\documentclass[onecolumn,twoside]{IEEEtran}
\usepackage{mathpazo}
\usepackage{times}
\usepackage{amsmath}
\usepackage{amsfonts}
\usepackage{latexsym}
\usepackage{amssymb}
\usepackage{upref}
\usepackage{theorem}
\usepackage{graphicx}
\usepackage{psfrag}
\usepackage{cite}
\usepackage{enumerate}
\usepackage{color}

\usepackage[normalem]{ulem}

\usepackage{yhmath}

\usepackage{tikz}
\usepackage{caption}
\usepackage{subcaption}
\usepackage{overpic}

\theoremstyle{plain}
\theorembodyfont{\normalfont\slshape}

\newtheorem{thm}{Theorem$\!$}
\newenvironment{theorem}
{\begin{thm}\hspace*{-1ex}{\bf.}}{\end{thm}}

\newtheorem{clm}[thm]{Claim$\!$}

\newtheorem{lem}[thm]{Lemma$\!$}
\newenvironment{lemma}{\begin{lem}\hspace*{-1ex}{\bf.}}{\end{lem}}

\newtheorem{prop}[thm]{Proposition$\!$}
\newenvironment{proposition}{\begin{prop}\hspace*{-1ex}{\bf.}}{\end{prop}}

\newtheorem{cor}[thm]{Corollary$\!$}

\newtheorem{defn}[thm]{Definition$\!$}
\newenvironment{definition}{\begin{defn}\hspace*{-1ex}{\bf.}}{\hfill $\Box$ \end{defn}}

\newtheorem{xmpl}[thm]{Example$\!$}
\newenvironment{example}{\begin{xmpl}\hspace*{-1ex}{\bf.}}{\hfill $\Box$ \end{xmpl}}

\newtheorem{cnstr}{Construction$\!$}

\newtheorem{rmk}[thm]{Remark$\!$}
\newenvironment{remark}{\begin{rmk}\hspace*{-1ex}{\bf.}}{\hfill $\Box$ \end{rmk}}

\newcounter{enumrom}
\renewcommand{\theenumrom}{(\roman{enumrom})}

\makeatletter
\renewcommand{\@endtheorem}{\endtrivlist}
\makeatother

\makeatletter
\renewcommand{\thefigure}{{\@arabic\c@figure}}
\renewcommand{\fnum@figure}{{\bf Figure\,\thefigure}}
\makeatother

\newcommand{\cB}{\mathcal{B}}

\newcommand{\cF}{\mathcal{F}}

\newcommand{\cL}{\mathcal{L}}

\newcommand{\cP}{\mathcal{P}}

\newcommand{\cX}{\mathcal{X}}

\newcommand{\mathset}[1]{\left\{#1\right\}}
\newcommand{\abs}[1]{\left|#1\right|}
\newcommand{\dabs}[1]{\left\|#1\right\|}

\newcommand{\floorenv}[1]{\left\lfloor #1 \right\rfloor}

\newcommand{\parenv}[1]{\left( #1 \right)}
\newcommand{\sparenv}[1]{\left[ #1 \right]}

\newcommand{\ip}[1]{\left\langle #1 \right\rangle}

\newcommand{\be}[1]{\begin{equation}\label{#1}}
\newcommand{\ee}{\end{equation}}
\newcommand{\eq}[1]{(\ref{#1})}

\renewcommand{\le}{\leqslant}
\renewcommand{\leq}{\leqslant}
\renewcommand{\ge}{\geqslant}
\renewcommand{\geq}{\geqslant}

\renewcommand{\Bbb}{\mathbb}

\newcommand{\Tref}[1]{Theo\-rem\,\ref{#1}}
\newcommand{\Pref}[1]{Pro\-po\-si\-tion\,\ref{#1}}
\newcommand{\Lref}[1]{Lem\-ma\,\ref{#1}}
\newcommand{\Cref}[1]{Co\-ro\-lla\-ry\,\ref{#1}}

\renewcommand{\Bbb}{\mathbb}

\newcommand{\C}{{\Bbb C}}
\newcommand{\N}{{\Bbb N}}

\newcommand{\R}{{\Bbb R}}
\newcommand{\Z}{{\Bbb Z}}
\newcommand{\bP}{\mathbb{P}}
\newcommand{\E}{\mathbb{E}}

\newcommand{\cPsi}{\cP_{\mathrm{si}}}
\newcommand{\ocP}{\overline{\cP}}
\newcommand{\tcP}{\overline{\cP}}
\newcommand{\God}{\Gamma^{\otimes d}}
\newcommand{\oGod}{\Gamma^{\boxtimes d}}
\newcommand{\oGodk}{\Gamma_{k,p}^{\boxtimes d}}
\newcommand{\hind}{\widehat{h_{\mathrm{ind}}}}
\newcommand{\ohind}{h_{\mathrm{ind}}}
\newcommand{\hindcom}{\ohind^{\mathrm{com}}}
\newcommand{\thind}{h_{\mathrm{ind}}}

\DeclareMathOperator{\fr}{fr}

\DeclareMathOperator{\sint}{int}

\newcommand{\shift}[1]{\sigma_{\mathbf{#1}}}
\newcommand{\tcap}{\mathsf{cap}}
\newcommand{\ccap}{\widehat{\tcap}}
\newcommand{\capcom}{\tcap^{\mathrm{com}}}
\newcommand{\cBcom}{\cB^{\mathrm{com}}}

\newcommand{\limsupup}[1]{\limsup_{#1\rightarrow\infty}}

\newcommand{\eqdef}{\triangleq}
\newcommand{\ZD}{{\Z^d}}
\newcommand{\Be}{\mathbb{B}_{\epsilon}}
\newcommand{\Bd}{\mathbb{B}_{\delta}}
\newcommand{\otid}{\otimes_{i\in [d]}}
\newcommand{\hpi}{\overline{\pi}}
\newcommand{\hfr}{\hat{\fr}}
\newcommand{\bu}{\mathbf{u}}
\newcommand{\bv}{\mathbf{v}}

\newcommand{\bunt}{\mathbf{e}}
\newcommand{\bone}{\mathbf{1}}
\newcommand{\bzer}{\mathbf{0}}
\newcommand{\given}{~\middle|~}

\outer\def\proclaim #1. #2\par{\medbreak
 \noindent{\bf#1.\enspace}{\sl#2\par}%
 \ifdim\lastskip<\medskipamount \removelastskip\penalty55\medskip\fi}

\begin{document}

\title{\textbf{On Independence And Capacity of Multidimensional Semiconstrained Systems}}

\author{\large
Ohad Elishco,~\IEEEmembership{Student Member,~IEEE},
Tom Meyerovitch,
Moshe~Schwartz,~\IEEEmembership{Senior Member,~IEEE}
\thanks{
  The material in this paper was presented in part at the
  IEEE International Symposium on Information
  Theory (ISIT 2017), Aachen, Germany, June 2017.}%
\thanks{Ohad Elishco is with the Department
  of Electrical and Computer Engineering, Ben-Gurion University of the Negev, Beer Sheva 8410501, Israel
   (e-mail: ohadeli@post.bgu.ac.il).}%
\thanks{Tom Meyerovitch is with the Department of Mathematics, Ben-Gurion University of the Negev,
   Beer Sheva 8410501, Israel (e-mail: mtom@math.bgu.ac.il).}%
\thanks{Moshe Schwartz is with the Department
   of Electrical and Computer Engineering, Ben-Gurion University of the Negev,
   Beer Sheva 8410501, Israel
   (e-mail: schwartz@ee.bgu.ac.il).}%
\thanks{This work was supported in part by the People Programme (Marie Curie Actions) of the European Union's Seventh Framework Programme (FP7/2007-2013) under REA grant agreement no.~333598 and by the Israel Science Foundation (grant no.~626/14).} 
}

\maketitle

\begin{abstract}
  We find a new formula for the limit of the capacity of
  certain sequences of multidimensional semiconstrained systems as
  the dimension tends to infinity. We do so by generalizing the
  notion of independence entropy, originally studied in the context of
  constrained systems, to the study of semiconstrained systems. Using
  the independence entropy, we obtain new lower bounds on the capacity
  of multidimensional semiconstrained systems in general, and
  $d$-dimensional axial-product systems in particular. In the case of
  the latter, we prove our bound is asymptotically tight, giving the
  exact limiting capacity in terms of the independence entropy. We
  show the new bound improves upon the best-known bound in a case
  study of $(0,k,p)$-RLL.
\end{abstract}

\section{Introduction}
\label{sec:intro}
\IEEEPARstart{E}{rror}-correcting codes and constrained codes may be
considered as two extreme ways of coping with a noisy channel. The
former are usually data independent, and assume errors are a
statistical phenomenon, reducing data-transmission rate to protect
against such errors. Constrained codes, however, assume certain
patterns in the data stream are responsible for the occurrence of
errors. Thus, constrained codes eliminate all undesirable patterns, at
the cost of reduced data-transmission rate.

Recently in \cite{EliMeySch16,EliMeySch16a}, semiconstrained systems
(SCSs) were suggested as a generalization to constrained systems
(which we emphasize by calling fully constrained systems). In SCSs we
do not eliminate the undesirable patterns entirely but rather we allow
them to appear with a restriction on their
frequency. To illustrate, consider a
  binary channel in which the appearance of $k$-consecutive $1$'s is
  forbidden. The set of allowed words is the well known inverted
  $(0,k)$-RLL. However, if $k$-consecutive $1$'s are not forbidden
  entirely, but instead are allowed to appear in at most a fraction
  $p$ of places, then the set of allowed words forms a SCS called the
  $(0,k,p)$-RLL system. Informally, a SCS is defined by a set
$\Gamma$ of probability measures over $k$-tuples. The allowed words in
the SCS are those in which the empirical distribution of $k$-tuples
belongs to $\Gamma$. This may be viewed as a generalization of fully
constrained systems since taking $\Gamma$ to be a subset with a
$0$-frequency restriction on some $k$-tuples yields a fully
constrained system.

SCSs not only generalize fully constrained systems, but also subsume a
range of other settings, which were mainly dealt with in an ad-hoc
fashion. Among these we can find DC-free RLL coding \cite{Kur11},
constant-weight ICI coding for flash memories
\cite{KaySie14,QinYaaSie14,CheChrKiaLinNguVu16a,CheChrKiaLinNguVu16b},
coding to mitigate the appearance of ghost pulses in optical
communication \cite{ShaSkiTur10,ShaTurTur07}, and the more general,
channel with cost constraints \cite{KarNeuKha88,KhaNeu96}.

In the one-dimensional case, the capacity of a SCS is given by a
relatively explicit expression as the solution to a certain
optimization problem on a finite dimensional space, e.g.,
\cite{MarRot92}. A probabilistic encoder for SCSs was constructed in
\cite{EliMeySch16}, and constant-bit-rate to constant-bit-rate
encoders are possible by approximating a SCS with a fully constrained
system, as described in \cite{EliMeySch16a}.

A natural extension, and the goal of this work, is to study
multidimensional SCSs. This is an extremely challenging problem,
considering the fact that even for fully constrained systems in
complete generality it is provably impossible to find an exact
solution. The capacity of multidimensional fully constrained systems
is known exactly only in a handful of cases
\cite{Lie67,Bax80,SchBru08,LouMar10}. In the absence of a general
method for computing the capacity, various bounds and approximations
were studied, e.g.,
\cite{CalWil98,SieWol98,ItoKatNagZeg99,KatZeg99,NagZeg00,OrdRot04,HalCheRotSieWol04,ShaRot10a,TalRot10,SchVar11}. It
should be emphasized that apart from its independent intellectual
merit, studying multidimensional systems is of practical importance
since most storage media are two- or three-dimensional, including
magnetic recording devices such as hard drives, optical recording
devices such as CDs and DVDs, and flash memories.

The approach we take in this work is bounding the capacity by studying
the independence entropy of SCSs, thus extending the works
\cite{LouMarPav13,MeyPav14}. The independence entropy appeared in
previous works on $d$-dimensional shifts of finite type.  Although
this notion was first defined in \cite{LouMarPav13}, the idea stemmed
from tradeoff functions studied in \cite{PooChaMar06}. It was defined
in a combinatorial fashion, where in this work we redefine it in a
probabilistic fashion. We show that the two definitions are equal for
the special case of fully constrained systems.
	
The motivation for the use of independence entropy is the fact that it
is more easily computable, since we only need to consider independent
probability measures which satisfy the constraints. We also focus on
the class of $d$-dimensional axial-product constraints, which form
a significant proportion of multidimensional fully
constrained systems studied thus far. For this class, our approach has
an additional major advantage in that instead of calculating the
independence entropy for a $d$-dimensional axial product SCS, we may
calculate it directly from the one-dimensional system. This
dimensionality reduction offers further simplification of the
calculations.

There are new features and difficulties that come up when
  adapting the results from fully constrained systems.  In an abstract
  sense, a very useful property of fully constrained systems is the
  following: If a measure $\mu$ is contained in some fully constrained
  system, and $\mu$ is a convex combination of measures, then each of
  them is contained in the same fully constrained system. This
  property does not hold for general semiconstrained systems. This is
  manifested for instance in the fact that any subword of an
  admissible word in a fully constrained system is also admissible,
  leading to sub-additivity of the sequence of the amount of
  admissible words. This, in turn, allows the use of Fekete's Lemma.

The main contributions of this paper are a formulation of the
independence entropy for SCSs, and its study in relation to the
capacity of SCSs. As a result, we obtain a new lower bound on the 
capacity of multidimensional SCSs, generalizing the results
of \cite{LouMarPav13,MeyPav14}, and in an example test case, improving
upon the best known bounds on the capacity of multidimensional
$(0,1,p)$-RLL SCSs given in \cite{EliMeySch16}.

In this work we also establish an equality of the limiting capacity
and independence entropy for the $d$-axial-product SCSs. As the
independence entropy is a lower bound on the entropy of a given SCS in
every dimension, the capacity approaches the independence
entropy as the dimension grows.

This paper is organized as follows. In Section \ref{sec:pre} we
describe the notation and give the required definitions used
throughout the paper. In Section \ref{sec:ind} we define the
independence entropy and provide results characterizing the
independence entropy. In Section \ref{sec:lower} we show that the 
capacity is lower bounded by the independence entropy. In
Section \ref{sec:upper} we show that the limiting capacity
of the $d$-axial-product SCS is equal to the independence entropy. We
conclude in Section \ref{sec:conc} by describing a short case study,
and comparing it with previous results. The appendices provide proofs
that the generalized notions we define in this paper indeed contain
fully constrained systems as a special case, thus providing a
generalization for them.

\section{Preliminaries}
\label{sec:pre}

Let $\N$ denote the set of natural numbers. We use $\bunt_i$ to denote
the unit vector of direction $i$, $\bzer$ to denote the all-zero
vector, and $\bone$ to the denote the all-one vector, where in all
cases, the dimension of the vectors is implied by the context. For
$n\in\N$ we define
\[[n]\eqdef\mathset{0,1,\dots,n-1}.\]
We shall often use $[n]\bunt_i$ to denote the set $\mathset{0\cdot
  \bunt_i,1\cdot \bunt_i,\dots,(n-1)\cdot \bunt_i}$. For $d,n\in\N$,
denote by $F_n^d$ the $d$-dimensional cube of length $n$, i.e., the
set $F_n^d\eqdef [n]^d$.  Obviously $|F_n^d|=n^d$. Additionally, for
$(n_0,\dots,n_{d-1})\in\N^d$ we conveniently denote
\[[(n_0,\dots,n_{d-1})]\eqdef [n_0]\times[n_1]\times\dots\times[n_{d-1}].\]

Throughout the paper, $\Sigma$ will be used to denote a finite
alphabet.  A word (or block) $w$ of length $n$ is a sequence of $n$
letters from $\Sigma$, denoted $w=a_0 a_1\dots a_{n-1}$, with
$a_i\in\Sigma$. We let $\abs{w}$ denote the length of the word $w$. We
can also consider infinite-sized words by mapping letters from
$\Sigma$ to positions on the integer grid $\ZD$. Such a word will be
denoted by $x\in\Sigma^\ZD$, and the letter in the $\bv\in\ZD$
position will be denoted by $x_{\bv}$ (sometimes referred to as the
restriction of $x$ to $\bv$). More generally, given any subset of the
integer grid, $S\subseteq\ZD$, a word $x\in \Sigma^S$ is a mapping of
letters from $\Sigma$ to positions indexed by elements of $S$.

We require a notation for sets of probability measures and their
marginals. For a set $W$ we denote by $\cP(W)$ the set of all
probability measures over $W$.

\begin{definition}
  Let $(X,\cB)$ be a measurable space. For every $\mu,\nu\in\cP(X)$,
  the \emph{total variation distance} is defined as
  \[
  \dabs{\mu-\nu}_{TV}\eqdef\sup_{W\in\cB} \abs{\mu(W)-\nu(W)}.
  \]
\end{definition}

Given a compact topological space $X$, the space $\cP(X)$ is itself a
compact topological space with respect to the weak $*$-topology.  In
particular, when $X$ is a finite set with the discrete topology, the
topology on $\cP(X)$ is given by the total variation distance which
also satisfies $\dabs{\mu-\nu}_{TV}=\frac{1}{2}\sum_{x\in
  X}\abs{\mu(x)-\nu(x)}$.

Given a continuous map $f:X \to Y$ between topological spaces, and
$\mu \in \cP(X)$, let $f(\mu) \in \cP(Y)$ be given by
\[f (\mu)(W) \eqdef \mu(f^{-1}(W)), ~ W \subseteq Y.\]

\begin{definition}
  For $d\in\N$, $S \subseteq \tilde{S} \subseteq \ZD$, and $x \in
  \Sigma^{\tilde S}$, let $x_S$ denote the restriction of $x$ to the
  coordinates in $S$. Let $\pi_S^{\tilde{S}} : \Sigma^{\tilde{S}} \to
  \Sigma^{S}$ denote the restriction map given by
  \[\pi_S^{\tilde{S}}(x)\eqdef x_S.\]
  When $\tilde{S}$ is clear from the context, we shall write $\pi_S$
  instead of $\pi_S^{\tilde{S}}$.
\end{definition}

While having the notation $\pi_S(x)$ in addition to the equivalent
notation $x_S$, seems superfluous, we shall require the former to
simplify our presentation. As a consequence of the previous
definition, for $\mu \in \cP(\Sigma^{\tilde{S}})$ and $S \subseteq
\tilde{S}$, we note that $\pi_{S}(\mu) \in \cP(\Sigma^{S})$ is the
$S$-marginal of $\mu$.

\begin{definition}
For $d\in\N$, $\bv \in \ZD$, let $\sigma_{\bv}: \Sigma^{\ZD} \to
\Sigma^{\ZD}$ be the shift by the vector $\bv$, given by
\[  (\sigma_{\bv}(x))_{\bu} \eqdef  x_{\bu+\bv},~ \bu \in \ZD,~ x \in \Sigma^{\ZD}.\]
We denote by $\cPsi(\Sigma^{\ZD})$ the space of shift-invariant
probability measures on $\Sigma^{\ZD}$, namely,
\[
\cPsi(\Sigma^{\ZD}) \eqdef \mathset{\mu \in \cP(\Sigma^{\ZD}) ~:~ \text{$\sigma_\bv(\mu) =\mu$ for all $\bv \in \ZD$} }.
\]
For $k \in \N$ we say that $\mu \in \cP(\Sigma^{F_k^d})$ is shift
invariant if it is the projection of some shift-invariant measure on
$\Sigma^{\ZD}$, i.e., if there exists $\tilde{\mu} \in
\cPsi(\Sigma^{\ZD})$ such that $\mu = \pi_{F_k^d} \tilde{\mu}$.  We denote
by $\cPsi(\Sigma^{F_k^d})$ the space of shift-invariant probability
measures on $\Sigma^{F_k^d}$, namely,
\[\cPsi(\Sigma^{F_k^d})\eqdef \pi_{F_k^d}(\cPsi(\Sigma^{\ZD}))
\subseteq \cP(\Sigma^{F_k^d}).\]
\end{definition}

In the one-dimensional case, $d=1$, it is rather easy to check whether
a given probability measure $\mu \in \cP(\Sigma^{F_k^1})$ is shift
invariant. Indeed, $\mu \in \cPsi(\Sigma^{F_k^1})$ if and only if it
satisfies the following finite system of linear equations,
\[\sum_{a\in\Sigma}\mu(a,a_1,\dots,a_{k-1})=\sum_{a\in\Sigma}\mu(a_1,\dots,a_{k-1},a),\]
for all $a_1,\dots,a_{k-1}\in\Sigma$.  

When $d \geq 2$ the space of finite marginals of shift invariant measures becomes much more complicated. It is
still not difficult to formulate an analogous system of linear
equations that are satisfied for every $\mu \in
\cPsi(\Sigma^{F_k^d})$. However, these linear conditions are no longer
sufficient conditions for shift invariance. In fact, the problem of
checking whether a given $\mu \in \cP(\Sigma^{F_k^d})$ is shift
invariant, is undecidable (assuming some computable representation of $\mu$). See for instance \cite{ChaGamHocUga12}, and
references within, for a related discussion.

We are interested in defining empirical distributions of words. To
that end, we give some more general definitions that we then
specialize to our specific needs. Given $x \in \Sigma^{\ZD}$, the
delta measure at $x$, denoted by $\delta_x \in \cP(\Sigma^{\ZD})$, is
defined by $\delta_x(\mathset{x})=1$. Additionally, given $n\in\N$,
the \emph{empirical measure associated with $x$ and $n$}, denoted
$\fr_{x,n} \in \cP(\Sigma^{\ZD})$, is given by
\[\fr_{x,n}\eqdef \frac{1}{n^d}\sum_{\bv \in F_n^d}\delta_{\sigma_\bv(x)}.\]
For $S\subseteq \ZD$ we can take the $S$-marginal, and define $\fr_{x,n}^S
\in \cP(\Sigma^{S})$ by
\[\fr^S_{x,n} \eqdef \pi_S (\fr_{x,n}).\]
Any word $w \in \Sigma^{F_n^d}$ may be extended periodically to the entire integer grid $\hat{w}\in\Sigma^\ZD$ by defining
\[ \hat{w}_\bv \eqdef w_{\bv\bmod n}\]
for all $\bv\in\ZD$, and where the modulo is taken entry-wise. The
empirical distribution we shall be requiring may now be defined.

\begin{definition}
Let $d,n\in\N$, $w\in\Sigma^{F_n^d}$, and $S\subseteq\ZD$. The
\emph{empirical distribution} of $w$ with respect to $S$, denoted
$\fr^S_w$, is defined by
\[\fr^S_w \eqdef \fr^S_{\hat{w}, n}.\]
\end{definition}

Combinatorially speaking, the empirical distribution $\fr_w^S$ is
obtained by cyclically scanning $w$ with an $S$-shaped window and
recording the frequency of the $S$-tuples in $w$. Thus, for instance,
given a word $w=w_0\dots w_{n-1}\in\Sigma^n$, $w_i\in\Sigma$, and $a
\in \Sigma^k$ we have
\[\fr^{[k]}_w(a)=  \frac{1}{\abs{w}}\sum_{i=0}^{\abs{w}-1}\mathbb{1}_{a}(w_i\dots w_{i+k-1})\]
where all coordinate indices are taken modulo $\abs{w}$, and
$\mathbb{1}_{a}:\Sigma^k \to \mathset{0,1}$ is the indicator function
of the singleton $\mathset{a}$.

\begin{example}
  Let $\Sigma=\mathset{0,1}$ and let $w=0010111001\in
  \Sigma^{F_{10}^1}$. We have that $|F_{10}^1|=10$ and
  \begin{align*}
    \fr^{[3]}_w(110)&=\frac{1}{10}\sum_{i=0}^{9}\Bbb{1}_{110}(w_iw_{i+1}w_{i+2})=\frac{1}{10},\\
    \fr^{[2]}_w(10)&=\frac{1}{10}\sum_{i=0}^{9}\Bbb{1}_{10}(w_iw_{i+1})=\frac{3}{10}
  \end{align*}
\end{example}

\begin{example}
  Let $\Sigma=\mathset{0,1}$ and consider
  \[w=\begin{bmatrix} 0&1&1&1\\ 0&0&1&1\\1&0&0&1\\ 1&0&1&0 \end{bmatrix}\in \Sigma^{F_4^2},
  \qquad\qquad a=\begin{bmatrix}
  0&1\\ 1&0\end{bmatrix}\in \Sigma^{F_2^2}.\]
  Then $\fr_w^{F_2^2}(a)=\frac{2}{16}$
  since, of the sixteen $2\times 2$ windows, exactly two contain $a$,
  shown in bold in the following:
  \[
  \begin{bmatrix} 0&1&1&1\\ 0&0&1&1\\1&0&\bzer&\bone\\ 1&0&\bone&\bzer \end{bmatrix},
  \qquad\qquad
  \begin{bmatrix} \bzer&1&1&\bone\\ 0&0&1&1\\1&0&0&1\\ \bone&0&1&\bzer \end{bmatrix}.
  \]
\end{example}

\begin{lemma}
  \label{lem:wlcontr}
  Suppose $d,n \in \N$, $w \in \Sigma^{F_n^d}$, and $S \subseteq
  \tilde{S} \subseteq \ZD$. Then
  \[\pi^{\tilde{S}}_S(\fr^{\tilde{S}}_w)=\fr^S_w.\]
\end{lemma}

\begin{IEEEproof}
  Let us denote $\mu\eqdef \fr_{\hat{w},n}\in\cP(\Sigma^\ZD)$.  By
  definition, for the right-hand side of the claim, for every
  $W\subseteq \Sigma^S$,
  \[ \fr^S_w(W) = \pi^{\ZD}_S(\mu)(W) = \mu\parenv{(\pi^{\ZD}_S)^{-1}(W)}.\]
  Similarly, for the left-hand side,
  \[ \pi^{\tilde{S}}_S(\fr^{\tilde{S}}_w)(W) = \pi^{\tilde{S}}_S\parenv{ \pi^{\ZD}_{\tilde{S}} (\mu)}(W) = \pi^{\ZD}_{\tilde{S}} (\mu) \parenv{ (\pi^{\tilde{S}}_S)^{-1}(W)} = \mu\parenv{(\pi^{\ZD}_{\tilde{S}})^{-1}\parenv{ (\pi^{\tilde{S}}_S)^{-1}(W)}}.\]
  But clearly for all $A\subseteq\Sigma^S$,
  \[(\pi^{\ZD}_S)^{-1}(W) = (\pi^{\ZD}_{\tilde{S}})^{-1}\parenv{ (\pi^{\tilde{S}}_S)^{-1}(W)}.\]
\end{IEEEproof}

Lemma \ref{lem:wlcontr} implies that the empirical
frequency of $S$-tuples in $w$ can be calculated by first calculating
the empirical frequency of $\tilde{S}$-tuples, and then taking the
$S$-marginal.

\begin{example}
    Let $\Sigma=\mathset{0,1}$ and consider
  \[w=\begin{bmatrix} 0&1&1&1\\ 0&0&1&1\\1&0&0&1\\ 1&0&1&0 \end{bmatrix}\in \Sigma^{F_4^2}.\]
  Take $S=[1]^2=\mathset{(0,0)}$ and
  $\tilde{S}=[(2,1)]=\mathset{(0,0),(1,0)}$. Then
  \begin{align*}
    \fr^{\tilde{S}}_w(00)=\frac{2}{16},\qquad \fr^{\tilde{S}}_w(01)=\frac{5}{16},\\
    \fr^{\tilde{S}}_w(10)=\frac{5}{16},\qquad \fr^{\tilde{S}}_w(11)=\frac{4}{16}.
  \end{align*}
  Moreover, we have that 
  \begin{align*}
    \fr_w^{S}(0)=\frac{7}{16},\qquad \fr_w^S(1)=\frac{9}{16}.
  \end{align*}
  We can verify now that
  \[\pi^{\tilde{S}}_S(\fr_w^{\tilde{S}})(0)=\fr_w^{\tilde{S}}\parenv{(\pi^{\tilde{S}}_S)^{-1}(0)}=\fr_w^{\tilde{S}}\parenv{\mathset{00,01}}=\frac{7}{16}=\fr_w^S(0).\]
\end{example}

We are now ready to define multidimensional semiconstrained systems.
\begin{definition}
  \label{def:scs}
  For $d\in\N$, a $\ZD$-\emph{semiconstrained system} (SCS) is a set
  $\Gamma\subseteq\cP(\Sigma^S)$ for some finite set $S \subseteq
  \ZD$. For $n \in \N$, the \emph{admissible $n$-blocks} of $\Gamma$
  are
  \[\cB_n(\Gamma) \eqdef \mathset{w\in\Sigma^{F_n^d} ~:~ \fr^S_w\in\Gamma}.\]
\end{definition}

Since all SCSs we study in this paper are $\ZD$-SCSs, we shall
abbreviate and call them just SCSs, where the dimension, $d$, will be
clear from the context.

Note that SCSs generalize $d$-dimensional fully constrained
systems. Recall that fully constrained systems are defined by a set of
``forbidden patterns'', $A\subseteq \Sigma^{F_k^d}$, such that a word
$w\in\Sigma^\ZD$ is admissible if and only if none of the elements of
$A$ appear as an $F_k^d$-tuple of $w$. Thus, fully constrained systems
correspond to subshifts of finite type in symbolic dynamics. In our
notation, we therefore have the following.

\begin{definition}
  \label{def:fullynew}
  For $d,k\in\N$, we say that $\Gamma \subseteq \cP(\Sigma^{F_k^d})$
  is \emph{fully constrained} if there exists some $L \subseteq
  \Sigma^{F_k^d}$ such that
  \[\Gamma = \{ \mu \in \cP(\Sigma^{F_k^d}) ~:~ \mu(L) = 1\}.\]
\end{definition}

\begin{example}
  Let $\Sigma=\mathset{0,1}$, take
  \[ L = \Sigma^{F_2^2} \setminus \mathset{
    \begin{bmatrix} 0&0\\1&1 \end{bmatrix},
    \begin{bmatrix} 0&1\\1&1 \end{bmatrix},
    \begin{bmatrix} 1&0\\1&1 \end{bmatrix},
    \begin{bmatrix} 1&1\\1&1 \end{bmatrix},
    \begin{bmatrix} 1&0\\1&0 \end{bmatrix},
    \begin{bmatrix} 1&1\\1&0 \end{bmatrix}
  },\]
  and consider the fully constrained system, $\Gamma$, defined by
  \[\Gamma=\mathset{\mu\in\cP(\Sigma^{F_2^2}) ~:~ \mu(L)=1}.\]
  Note that $\cB_n(\Gamma)$ is the set of all $n\times n$
  two-dimensional binary arrays such that none of the six patterns
  above appears within a $2\times 2$ window in them. It is simple to
  verify that in fact, no two horizontally adjacent $1$'s may appear,
  and no two vertically adjacent $1$'s may appear, in any admissible
  word. Thus, the $n\times n$ arrays in $\cB_n(\Gamma)$ are the
  admissible words of the (cyclical) $(1,\infty)$-RLL fully
  constrained system.
\end{example}

An important figure of merit we associate with any set of words, and
in particular, with SCSs, is the capacity, which we now define.
\begin{definition}
    Let $d\in\N$, and let $S\subseteq\Z^d$ be a finite subset. For any SCS,
    $\Gamma \subseteq\cP(\Sigma^S)$, and for $\epsilon>0$, let
    \[\Be(\Gamma)\eqdef\mathset{\mu\in\cP(\Sigma^S) ~:~ \inf_{\nu\in \Gamma} \dabs{\mu-\nu}_{TV}\leq\epsilon }.\]
	The \emph{capacity} of $\Gamma$ is defined as,
	\[\tcap(\Gamma) \eqdef \lim_{\epsilon \to 0^+}\limsupup{n}\frac{1}{n^d}\log_2\parenv{\abs{\cB_n(\Gamma)}}.\]
  \end{definition}

  First, we mention that $\lim_{\epsilon\to 0^+}$ in the definition of the
capacity exists due to monotonicity, since $\abs{\cB_n(\Be(\Gamma))}$ is
non-increasing in $\epsilon$.

To avoid certain pathological scenarios,
  \cite{EliMeySch16,EliMeySch16a} defined sets of weakly-admissible
  words and their capacity. We contend that the capacity definition
  provided here is the proper multidimensional generalization of these
  definitions. Intuitively, the capacity measures the exponential
  growth rate of the number of words that ``almost'' satisfy the
  semiconstraints given by $\Gamma$. Additionally, it has the nice property
  that the capacity of a set $\Gamma$ is equal to the capacity of the
  closure of $\Gamma$.

At first glance this definition of capacity may seem odd. A naive definition, which we call the \emph{internal capacity}, might be as follows.

\begin{definition}
	Let $d\in\N$, $S\subseteq\Z^d$ finite, and $\Gamma \subseteq\cP(\Sigma^S)$ be a SCS. 
	The \emph{internal capacity} of $\Gamma$ is defined as  
	\[\ccap(\Gamma)\eqdef\limsupup{n}\frac{1}{n^d}\log_2\parenv{\abs{\cB_n(\Gamma)}}.\] 
\end{definition}

By definition we have  
\[ \tcap(\Gamma)=\lim_{\epsilon \to 0^+}\ccap(\Be(\Gamma))\]
which means that 
\begin{equation}
  \label{eq:ccaptcap}
  \ccap(\Gamma)\leq \tcap(\Gamma).
\end{equation}
We also observe that for some ``nice'' SCSs $\Gamma$, $\ccap(\Gamma)=\tcap(\Gamma)$. For instance, we have the following result for
one-dimensional SCSs.
\begin{theorem}\cite[Section 2]{EliMeySch16a}
  \label{oldtheo1}
  Let $k\in\N$, and $\Gamma\subseteq\cP(\Sigma^k)$ be convex and equal
  to the closure of its relative interior in $\cPsi(\Sigma^k)$. Then
  \[\tcap(\Gamma) = \ccap(\Gamma)=\log_2 |\Sigma|-\inf_{\eta\in\Gamma\cap\cPsi(\Sigma^k)} H(\eta|\mu)\] 
  where $H(\cdot|\cdot)$ is the relative entropy function, and $\mu$
  is defined by $\mu(\phi a)\eqdef
  \frac{1}{|\Sigma|}\sum_{a'\in\Sigma}\eta(\phi a')$ for all
  $\phi\in\Sigma^{k-1}$ and $a\in\Sigma$.
\end{theorem}

\begin{remark}
Consider the (compact) space $M=\cP(\Sigma^S)$ and let $C(M)$ be the
set of all closed (hence, compact) subsets of $M$. Thus, $C(M)$ is a
compact topological space (under the Hausdorff metric). Since
$\ccap(\Gamma)$ is monotone, the set of $\Gamma$s for which
$\ccap(\Gamma)\neq\tcap(\Gamma)$ is meager. In other words, if we
consider $\ccap(\Be(\Bd(\Gamma)))$ as a function of $\epsilon$,
$f(\epsilon)\eqdef\ccap(\Be(\Gamma))$, then
$\tcap(\Bd(\Gamma))=\ccap(\Bd(\Gamma))$ whenever $f$ is continuous in
$\delta$. Since $f$ is a monotone function, it is discontinuous on a
countable number of places. In practice, it means that if for a
specific $\Gamma$, $\tcap(\Gamma)\neq\ccap(\Gamma)$ an arbitrary small
change in $\Gamma$ will give an equality.
\end{remark}

\begin{remark}
  For a fully constrained system, $\Gamma\subseteq\cP(\Sigma^{F_k^d})$,
  non-emptyness of $\cB_n(\Gamma)$ for all $n>0$ is equivalent to the
  fact that the subshift of finite type
  \[\mathset{w\in\Sigma^{\ZD} ~:~ \forall \bv\in\ZD,
    \parenv{\sigma_{\bv}(w)}_{F_k^d}\in L},\] is not empty. Berger's
  Theorem \cite{Ber66} implies that it is undecidable whether a
  subshift of finite type is empty given $L$. Because (under reasonable assumptions on the representation) it is undecidable if a given multidimensional SCS is non-empty, it is difficult to understand what a SCS really looks like.
\end{remark}

At this point we pause to ponder the following: Note that the
definition of empirical frequency is cyclic (in the sense that
coordinates are taken modulo $n$) while in traditional fully
constrained systems it is not. This seems at odds with our claim of
SCSs generalizing fully constrained systems. The necessity of the
modulo in the definition of SCSs stems from working with the space of
shift-invariant measures and their associated admissible
words. Shift-invariant measures are defined over $\Z^d$, hence, it is
necessary to complete a word $w\in\Sigma^{F_n^d}$ to a word from
$\Sigma^{\Z^d}$.  We choose to do this completion periodically using
the modulo notion, extending $w$ to $\hat{w}$. This choice simplifies
the analysis which follows. We contend that with respect to this
issue, the capacity is more natural than the internal capacity, since it is equal to the
non-cyclic capacity of fully constrained systems. To avoid a lengthy
detour, the full details are provided in Appendix \ref{appA}.

Finally, we raise the question: what multidimensional SCSs are of
interest? If we examine the extensive literature for fully constrained
systems, a significant proportion of multidimensional fully constrained
systems are defined as an axial product of one-dimensional fully
constrained systems. Intuitively speaking, if we have a set of
``forbidden patterns'' defining a one-dimensional fully constrained
system, we can define its $d$-dimensional axial product by forbidding
these patterns along each dimension. We now formally define this for
the case of $d$-dimensional SCSs with slightly more generality. This definition generalizes the $d$-dimensional axial product defined in \cite{LouMarPav13}.

\begin{definition}
  Consider $S_0,\dots,S_{d-1} \subseteq \N$, with $0 \in S_i$ for all
  $i\in [d]$, and SCSs $\Gamma_i\subseteq\cP(\Sigma^{S_i})$. Denote
  $S\eqdef \bigcup_{i\in [d]} S_i\bunt_i \subseteq \ZD$. The
  \emph{$d$-axial-product SCS}, denoted $\otid\Gamma_i$, is defined by
  \[ \otid\Gamma_i \eqdef \mathset{\mu\in\cP(\Sigma^{S}) ~:~ \forall i\in[d],\; \pi_{S_i\bunt_i}(\mu) \in\Gamma_i}.\]
\end{definition}

It follows from the above definition, that for every $n\in\N$ we have
\[ \cB_n\parenv{\otid\Gamma_i}= \mathset{w\in F_n^d ~:~ \forall i\in[d],\; \fr^{S_i\bunt_i}_w\in\Gamma_i},\]
with coordinates taken modulo $n$.  Intuitively, the arrays of a
$d$-axial-product SCS satisfy that along the $i$th direction, the
empirical distribution of $S_i$-tuples is in $\Gamma_i$. Note that
$\otid\Gamma_i$ induces a set of measures over $\Sigma^{F_k^d}$ where
$k=\max_{i} \mathset{k_i ~:~ k_i\in S_i}$. Hence, we sometimes
consider a $d$-axial-product SCS $\otid\Gamma_i$ as a subset of
$\cP(\Sigma^{F_k^d})$.

\begin{example}
  Let $\Sigma=\mathset{0,1}$. Consider two real constants $0\leq
  p_0,p_1\leq 1$, and the one-dimensional SCSs, $\Gamma_0$ and
  $\Gamma_1$, given by
  \begin{align*}
    \Gamma_0&=\mathset{\mu\in\cP(\Sigma^2) ~:~ \mu(11)\leq p_0}, \\
    \Gamma_1&=\mathset{\mu\in\cP(\Sigma^2) ~:~ \mu(11)\leq p_1}.
  \end{align*}
  Here we are taking $S_0=S_1=\mathset{0,1}$. The admissible words in
  the $2$-axial-product SCS, $\Gamma_0\otimes\Gamma_1$, are all
  two-dimensional words in which the empirical frequency of two
  horizontally adjacent $1$s is at most $p_0$, and the empirical
  frequency of two vertically adjacent $1$s is at most $p_1$, i.e.,
  all the words $w\in\Sigma^{F_n^2}$ such that
  \begin{align*}
  \fr_w^{\mathset{(0,0),(1,0)}}(11) &\leq p_0 ,\\
  \fr_w^{\mathset{(0,0),(0,1)}}(11) &\leq p_1.
  \end{align*}
  We may also consider $\Gamma_0\otimes\Gamma_1$ as a subset of
  $\cP(\Sigma^{F_2^2})$
  \[\Gamma_0\otimes\Gamma_1=\mathset{\mu\in\cP(\Sigma^{F_2^2}) ~:~ \pi_{\mathset{(0,0),(0,1)}}(\mu)(11)\leq p_0,\; \pi_{\mathset{(0,0),(1,0)}}(\mu)(11)\leq p_1}.\]
  Note that
  \begin{align*}
    \pi_{\mathset{(0,0),(1,0)}}(\mu)(11) &=
    \mu\parenv{\begin{bmatrix}0&0\\ 1&1\end{bmatrix}}
    +\mu\parenv{\begin{bmatrix}0&1\\ 1&1\end{bmatrix}}
    +\mu\parenv{\begin{bmatrix}1&0\\ 1&1\end{bmatrix}}
    +\mu\parenv{\begin{bmatrix}1&1\\ 1&1\end{bmatrix}},\\
    \pi_{\mathset{(0,0),(0,1)}}(\mu)(11)&=\mu\parenv{\begin{bmatrix}1&0\\ 1&0\end{bmatrix}}
    +\mu\parenv{\begin{bmatrix}1&0\\ 1&1\end{bmatrix}}
    +\mu\parenv{\begin{bmatrix}1&1\\ 1&0\end{bmatrix}}
    +\mu\parenv{\begin{bmatrix}1&1\\ 1&1\end{bmatrix}}.
  \end{align*}
\end{example}

In this paper we are interested in the capacity and the internal
capacity of multidimensional SCSs. Although the capacity is
easier to work with, as we will see later on, the task of computing it
is still daunting. Thus, there is a necessity for more easily
computable bounds on the capacity. To this end, we define the
independence entropy of a $d$-dimensional SCS, which is the basis of
the main results of this paper.

\section{Independence Entropy}
\label{sec:ind}
In this section we define the independence entropy of
multidimensional SCSs and present some of its properties. It will be
used to bound the capacity. The independence entropy is not a
new notion, and has appeared previously in \cite{LouMarPav13} in
relation to the capacity of fully constrained systems. However, the
formulation of the independence entropy was combinatorial and
therefore less suitable for our purposes. Thus, we modify the definition
of independence entropy and formulate it as a statistical notion.

The admissible words of SCSs (see Definition \ref{def:scs}) have their
empirical $S$-tuple distribution from $\Gamma$. Finding such words
inexorably involves intricate dependencies between coordinates. This
affects not only the task of generating such words, but also the very
basic problem of calculating or bounding the capacity of the SCS --
the problem that is the focus of this paper.

In an attempt to simplify this problem, we study the
independence-entropy approach. We eliminate all dependencies by
considering only product measures, i.e., where the symbol in each
coordinate of the word is chosen independently of other
coordinates. Accordingly, we only require the \emph{average} of
$S$-marginals to be in $\Gamma$. We then ask what is the entropy of
such a system. Intuitively, we are seeking the maximum rate of
transmission in a system where word coordinates are transmitted
independently and in parallel, designed such that the average
$S$-marginals are in $\Gamma$. The following model provides a rough interpretation of the independence entropy: Suppose each bit of the output is transmitted by a different agent, and the number of agents is very large. The agents are allowed to coordinate a protocol in advance, but are unable to communicate once they receive the messages to be transmitted. In addition, the statistics of the output should roughly satisfy the constraints given by $\Gamma$, with high probability (as a function of the number of agents). In this case under suitable assumptions, the maximal transmission rate would coincide with the independence entropy.
We proceed with formal definitions,
starting with a product measure.

\begin{definition}
  Let $d\in\N$, and let $S \subseteq \ZD$ be a finite set. We say that
  $\mu\in\cP(\Sigma^{S})$ is \emph{an independent probability measure}
  or \emph{a product measure} if $\mu(w)=\prod_{\bv\in S}
  \pi_{\mathset{\bv}}(\mu)(w)$.  For $S \subseteq \ZD$ that is
  possibly infinite, $\mu\in\cP(\Sigma^{S})$ is a product measure
  whenever $\pi_{S'}(\mu)$ is a product measure for every finite $S'
  \subseteq S$.
\end{definition}

In other words, we say that $\mu$ is independent if there exists
$\mathset{p_{\bv}\in\cP(\Sigma) ~:~ \bv\in S}$ such that
$\mu=\prod_{\bv\in S}p_{\bv}$.  We naturally identify the set of
product measures in $\cP(\Sigma^S)$ with $\parenv{\cP(\Sigma)}^{S}$.

Next, we define the average of a marginal.
\begin{definition}
  Given $d,n\in\N$, $\mu \in \cP(\Sigma^{F_n^d})$, and $S \subseteq
  F_n^d$, let $\hpi_S(\mu) \in \cP(\Sigma^S)$ be the average of the
  $S$-marginals over translates of $\mu$:
  \[\hpi_S(\mu)\eqdef \frac{1}{|F_n^d|}\sum_{\bv\in F_n^d} \pi_{S+\bv} (\mu),\] 
  where the coordinates $S+\bv$ are taken modulo $n$.
\end{definition}

Let $S \subseteq F_k^d$ and let $\Gamma \subseteq \cP(\Sigma^{S})$ be
a SCS. For $n \geq k$ we define
\[
\ocP_n(\Gamma)\eqdef\mathset{
\mu\in\parenv{\cP(\Sigma)}^{F_n^d} ~:~
\hpi_{S}(\mu)\in\Gamma }.
\]
Thus, $\ocP_n(\Gamma)$ consists of product measures on
$\Sigma^{F_n^d}$ such that the average of the $S$-marginals is in
$\Gamma$.  We can now define the independence entropy of a SCS.

\begin{definition}
  Let $d,k\in\N$, $S\subseteq F_k^d$, and let
  $\Gamma\subseteq\cP(\Sigma^{S})$ be a $d$-dimensional SCS. The
  \emph{internal independence entropy} of $\Gamma$ is defined by
  \[\hind(\Gamma)\eqdef \limsupup{n} \sup_{\mu\in\ocP_n(\Gamma)} \frac{1}{n^d}H(\mu),\]
  where $H(\mu)\eqdef -\sum_{w\in\Sigma^{F_n^d}}\mu(w)\log_2 \mu(w)$ is the
  entropy of $\mu$. The \emph{independence entropy} of
  $\Gamma$ is defined by
  \[\ohind(\Gamma)\eqdef \lim_{\epsilon \to 0^+}\hind(\Be(\Gamma)).\]
\end{definition}

Again, it is clear by definition that
\begin{equation}
  \label{eq:hindohind}
  \hind(\Gamma)\leq \ohind(\Gamma).
\end{equation}
The notion of independence entropy which appears here is a
generalization of the combinatorial notion for fully constrained
systems that appears in \cite{LouMarPav13}.

\begin{theorem}
  \label{indentropyequivalent} 
  Let $d,k\in\N$, and let $\Gamma\subseteq\cP(\Sigma^{F_k^d})$ be a
  fully constrained system. Then
  \[\ohind(\Gamma)=\hindcom(\Gamma)\]
  where $\hindcom$ is the combinatorial independence
  entropy from \cite{LouMarPav13}.
\end{theorem}

To avoid a significant diversion from the main discussion, the proof
of \Tref{indentropyequivalent}, together with the required definitions
from \cite{LouMarPav13}, are given in Appendix \ref{appB}.

We now show properties of $\hind$ and $\ohind$ which make them easier
to analyze by reducing the multidimensional case to the
one-dimensional case. We start with an inequality given in the following lemma. The proof follows the same argument that was used in \cite{LouMarPav13} to show the inequality for fully constrained systems. However, the equality for fully constrained systems holds in an easier and stronger sense.

\begin{lemma}\label{lem:ind_axial}
  Let $k\in\N$, and let $\Gamma\subseteq\cP(\Sigma^k)$ be a
  one-dimensional SCS. Then for all $d\in\N$,
  \[\hind(\Gamma) \leq \hind(\Gamma^{\otimes d}).\]
\end{lemma}
\begin{IEEEproof}
  Take $\hat{\mu}\in\ocP_{n}(\Gamma)$. Since $\hat{\mu}$ is
  a product measure, it can be written as
  $\hat{\mu}=\prod_{i=0}^{n-1}\pi_{\mathset{i}}(\hat{\mu})$. We now
  construct a measure $\mu\in\ocP_n(\Gamma^{\otimes d})$
  using $\hat{\mu}$. For every $\bv\in F_n^d$ set
  \[\pi_{\{\bv\}}(\mu)\eqdef\pi_{\{\ell(\bv)\}}(\hat{\mu}),\]
  where $\ell(\bv) \eqdef\parenv{\sum_{i=0}^{d-1}v_i}\bmod n$ is the
  modulo $n$ of the sum of the coordinates of $\bv$.


  Observe that $\mu$ is such that in every row in every direction,
  i.e., a set of coordinates of the form $\bv+[n]\bunt_i$, we obtain
  some cyclic rotation of $\hat{\mu}$ by $t$ positions, denoted
  $\sigma_t(\hat{\mu})$. However,
  $\hat{\mu}\in\ocP_{n}(\Gamma)$ implies
  $\sigma_t(\hat{\mu})\in\ocP_{n}(\Gamma)$. Thus, we obtain
  that $\mu\in \ocP_n(\Gamma^{\otimes d})$ and
  \[ \frac{1}{n} H(\hat{\mu})=\frac{1}{n^d} H(\mu).\]
  Since we are taking the supremum over all measures $\hat{\mu}$, we
  have $\hind(\Gamma)\leq \hind(\Gamma^{\otimes d})$.
\end{IEEEproof}

\begin{theorem}
  \label{th:equ}
  Let $k\in\N$, and let $\Gamma\subseteq\cP(\Sigma^k)$ be a
  one-dimensional SCS. Then for all $d\in\N$,
  \[\ohind(\Gamma^{\otimes d})= \ohind(\Gamma).\]
\end{theorem}

\begin{IEEEproof}
  We first show that $\ohind(\Gamma^{\otimes d})\leq
  \ohind(\Gamma)$. Fix $\delta >0$ and take
  $\mu\in\ocP_n(\Bd(\Gamma^{\otimes d}))$. Recall that
  \[\ocP_n(\Bd(\Gamma^{\otimes d}))\subseteq \cP(\Sigma^{F_n^d}).\]
  Let $(\bv_i)_{i \in [n^{d-1}]}$ be an enumeration of
  $\mathset{0}\times F_n^{d-1}$, i.e.,
  \[\{\bv_0,\dots,\bv_{n^{d-1}-1}\} = \mathset{0}\times F_n^{d-1}.\]
  For $i \in [n^{d-1}]$, define $\mu_i \in \cP(\Sigma^{n})$ by $ \mu_i
  \eqdef \pi_{[n]\bunt_0+\bv_i}(\mu)$.  Now let
  $\hat{\mu}\in\cP(\Sigma^{n^d})$ be the product measure that is the
  product of all the $\mu_i$'s. This means that for a word
  $a=a_0\dots a_{n^d-1}\in\Sigma^{n^d}$,
  \[\hat{\mu}(a)\eqdef \mu_0(a_0\dots a_{n-1})\mu_1(a_n\dots a_{2n-1})\cdots\mu_{n^d-1}(a_{n(n^{d-1}-1)}\dots a_{n^d-1}),\]
  Since each of the $\mu_i$'s is already a product measure,
  $\hat{\mu}\in\cP(\Sigma^{n^d})$ is also a product measure. We have
  \begin{align*}
    \hpi_{[k]}(\hat{\mu})&=\frac{1}{n^d}\sum_{j=0}^{n^d-1}\pi_{j+[k]}(\hat{\mu})\\
    &=\frac{1}{n^d}\parenv{\sum_{i=0}^{n^{d-1}-1}\sum_{j=in}^{(i+1)n-1}\pi_{j+[k]}(\hat{\mu})}\\
    &=\frac{1}{n^d}\parenv{\sum_{i=0}^{n^{d-1}-1}\parenv{\sum_{j=in}^{(i+1)n-k}\pi_{j+[k]}(\hat{\mu})+\sum_{j=(i+1)n-k+1}^{(i+1)n-1}\pi_{j+[k]}(\hat{\mu})}}\\
    &\stackrel{(a)}{=}\frac{1}{n^d}\parenv{\sum_{i=0}^{n^{d-1}-1}\sum_{j=in}^{(i+1)n-k}\pi_{(j-in)+[k]}(\mu_i)+\sum_{i=0}^{n^{d-1}-1}\sum_{j=(i+1)n-k+1}^{(i+1)n-1}\pi_{j+[k]}(\hat{\mu})}\\
    &=\frac{1}{n^d}\parenv{\sum_{i=0}^{n^{d-1}-1}\sum_{j=in}^{(i+1)n-1}\pi_{(j-in)+[k]}(\mu_i)-\sum_{i=0}^{n^{d-1}-1}\sum_{j=(i+1)n-k+1}^{(i+1)n-1}\pi_{(j-in)+[k]}(\mu_i)+\sum_{i=0}^{n^{d-1}-1}\sum_{j=(i+1)n-k+1}^{(i+1)n-1}\pi_{j+[k]}(\hat{\mu})}\\
    &=\frac{1}{n^{d-1}}\parenv{\sum_{i=0}^{n^{d-1}-1}\hpi_{[k]}(\mu_i)-\frac{1}{n}\sum_{i=0}^{n^{d-1}-1}\sum_{j=(i+1)n-k+1}^{(i+1)n-1}\parenv{\pi_{(j-in)+[k]}(\mu_i)-\pi_{j+[k]}(\hat{\mu})}}\\
    &=\hpi_{[k]\bunt_0}(\mu)-\frac{1}{n^d}\sum_{i=0}^{n^{d-1}-1}\sum_{j=(i+1)n-k+1}^{(i+1)n-1}\parenv{\pi_{(j-in)+[k]}(\mu_i)-\pi_{j+[k]}(\hat{\mu})}
  \end{align*}
  where $(a)$ follows from the definition of $\hat{\mu}$ and since the
  coordinates are taken modulo $n$ when calculating
  $\pi_{[k]}(\mu_i)$.  Each
  $\parenv{\pi_{(j-in)+[k]}(\mu_i)-\pi_{j+[k]}(\hat{\mu})}$ is a
  signed measure of total variation norm at most $1$. Therefore,
  \[ \dabs{\hpi_{[k]}(\hat{\mu})-\hpi_{[k]\bunt_0}(\mu)}_{TV} \le \frac{k}{n}.\]
  This means that
  $\hpi_{[k]}(\hat{\mu})\in\mathbb{B}_{\frac{k}{n}+\delta}(\Gamma)$. We
  obtained that for every $\epsilon>\delta>0$, and every $\mu\in
  \ocP_n(\Bd(\Gamma^{\otimes d}))$, we can find $n_0\in\N$ such that
  for every $n>n_0$,
  $\hat{\mu}\in\ocP_{n^d}\parenv{\Be(\Gamma)}$. Since $\mu$ and
  $\hat{\mu}$ are both product measures we have
  \begin{align*}
    H(\mu)&=\sum_{\bv\in F_n^d}H(\pi_{\mathset{\bv}}(\mu))\\
    &=\sum_{i\in [n^d]}H(\pi_{\{i\}}(\hat{\mu}))\\
    &=H(\hat{\mu}).
  \end{align*}
  This implies that for every $\epsilon>\delta >0$,
  \[\limsupup{n}\sup_{\mu\in\ocP_n(\Bd(\Gamma^{\otimes d}))}\frac{1}{n^d}H(\mu)\leq \limsupup{n} \sup_{\mu\in\ocP_{n^d}(\Be(\Gamma))}\frac{1}{n^d}H(\mu).\]
  We therefore obtain $\ohind (\Gamma^{\otimes d})\leq \hind
  (\Be(\Gamma))$ for every $\epsilon>0$. Taking the limit as
  $\epsilon\to 0^+$, by the definition of $\ohind(\Gamma)$ we have
  \[\ohind(\Gamma^{\otimes d})\leq \ohind(\Gamma).\]
	
  We now show the other direction. By Lemma \ref{lem:ind_axial}, For
  every $\delta >0$ we have
  \[\hind(\Bd(\Gamma))\leq \hind(\Bd(\Gamma)^{\otimes d}).\]
  By monotonicity of $\hind$ it thus follows that for every $\delta
  >0$,
  \[\ohind(\Gamma)\leq \hind(\Bd(\Gamma)^{\otimes d}).\]
  Now observe that for every $\epsilon >0$ there exists $\delta >0$ so
  that
  \[\Bd(\Gamma)^{\otimes d} \subseteq \Be\parenv{\Gamma^{\otimes d}	}.\]
  It follows that for every $\epsilon >0$
  \[\ohind(\Gamma)\leq \hind\parenv{\Be\parenv{
  \Gamma^{\otimes d}}}.\]
  Thus, by taking the limit $\epsilon\to 0^+$,
  \[\ohind(\Gamma)\leq \ohind\parenv{
  \Gamma^{\otimes d}}.\]
\end{IEEEproof}


We conclude this section by noting that Lemma \ref{lem:ind_axial} and
Theorem \ref{th:equ} show that for $\Gamma\subseteq\cP(\Sigma^k)$,
\begin{equation}
  \label{eq:summ2}
  \hind(\Gamma)\leq\hind(\Gamma^{\otimes d})\leq\ohind(\Gamma^{\otimes d})=\ohind(\Gamma).
\end{equation}

\section{Independence Entropy Lower Bounds The Capacity}
\label{sec:lower}

This section and the next explore the relationship between the
independence entropy and the capacity. In this section we show that
the capacity of any $d$-dimensional SCS (not necessarily an axial
product) is lower bounded by the independence entropy.

Before proceeding we require a simple lemma.

\begin{lemma}
  \label{lem:contr}
  Let $d,n\in\N$, and $S\subseteq F_n^d$, then $\pi_S$ and $\hpi_S$
  are contractions with respect to the total-variation distance, i.e.,
  for all $\mu,\nu\in\cP(\Sigma^{F_n^d})$,
  \begin{align*}
    \dabs{\pi_S(\mu)-\pi_S(\nu)}_{TV} & \leq \dabs{\mu-\nu}_{TV}, \\
    \dabs{\hpi_S(\mu)-\hpi_S(\nu)}_{TV} & \leq \dabs{\mu-\nu}_{TV}.
  \end{align*}
\end{lemma}

\begin{IEEEproof}
  For every $W \subseteq \Sigma^{S}$ we have
  \[
  \abs{\pi_S(\mu)(W)-\pi_S(\nu)(W)} = \abs{\mu(\pi_S^{-1}(W))-\nu(\pi_S^{-1}(W))}
  \leq \sup_{A' \subset \Sigma^S} \abs{\mu(W')-\nu(W')} = \dabs{\mu-\nu}_{TV}.
  \]
  Hence the function $\pi_{S+\bv}$ is a contraction for every $\bv \in
  F_n^d$. Then $\hpi_S$, being an average of contractions, is itself a
  contraction.
\end{IEEEproof}

We are now ready to state and prove the main result of this section --
a lower bound on the capacity. The corresponding result for fully constrained systems was obtained in \cite{LouMarPav13}.

\begin{theorem} 
  \label{indentropy1}
  Let $d\in\N$, $S \subseteq \ZD$ be a finite set, and let
  $\Gamma\subseteq\cP(\Sigma^S)$ be a SCS.  Then $\ohind(\Gamma)\leq
  \tcap(\Gamma)$.
\end{theorem}

\begin{IEEEproof}
  Fix $\delta >0$, $n \in\N$ such that $S\subseteq F_n^d$, and let
  $\mu\in\ocP_n(\Bd(\Gamma))$. For $m\in\N$, we have a
  natural identification isomorphism $\Sigma^{F_{nm}^d} \cong
  (\Sigma^{F_n^d})^{F_m^d}$ that identifies $\bv\in F_{nm}^d$ with the
  unique pair $\mathbf{r}\in F_n^d$ and $\mathbf{q}\in F_m^d$ such
  that $\bv=n\mathbf{q}+\mathbf{r}$. Consider the product measure
  $\mu^m \in \cP(\Sigma^{F_n^d})^{F_m^d} \subseteq
  \cP(\Sigma^{F_{nm}^d})$ satisfying
  \[\mu^m(\mathset{x}) = \prod_{ \bv \in F_m^d} \mu(\pi_{F_n^d}(\sigma_{n\bv}(x))).\]
  Note that since $\mu$ is a product measure, $\mu^m$ is also a
  product measure.

  For a word $w\in\Sigma^{F_{nm}^d}$, denote by $\hfr^{F_n^d}_w$ the
  empirical distribution of \emph{non-overlapping} $F_n^d$-tuples,
  i.e.,
  \[\hfr^{F_n^d}_w\eqdef\frac{1}{|F_m^d|}\sum_{\bu\in
    F_m^d}\delta_{\pi_{F_n^d}(\sigma_{n\bu}(\hat{w}))}.\]
  Additionally, observe that
  \[ \frac{1}{|F_n^d|}\sum_{\bv\in F_n^d} \hfr^{F_n^d}_{\sigma_{\bv}(w)} =\fr_w^{F_n^d}.\]
  Also, because $\pi_S^{F_n^d}$ is an affine map, it follows that 
  \[ \frac{1}{|F_n^d|}\sum_{\bv\in F_n^d} \pi_S^{F_n^d}\hfr^{F_n^d}_{\sigma_{\bv}(w)} =\pi_S (\fr_w^{F_n^d}). \]
  By Lemma \ref{lem:wlcontr}, $ \pi_S (\fr_w^{F_n^d})=\fr_w^{S}$.
  
  Note that by the construction of $\mu^m$ we have $\hpi_S(\mu)=\hpi_S(\mu^m)$, and we obtain,
  \begin{align*}
    \dabs{\fr^S_w-\hpi_S(\mu^m)}_{TV}&= \dabs{\fr^S_w-\hpi_S(\mu)}_{TV} \\
    &=\dabs{\frac{1}{|F_{nm}^d|}\sum_{\bu\in F_{nm}^d}\pi_S(\delta_{\sigma_{\bu}(w)})-\frac{1}{|F_n^d|}\sum_{\bv\in F_n^d}\pi_{S+\bv}(\mu)}_{TV}\\
    &=\dabs{\frac{1}{|F_n^d||F_{m}^d|}\sum_{\bv\in F_n^d}\sum_{\bu\in F_{m}^d} \pi_S(\delta_{\sigma_{n\bu+\bv}(w)})-\frac{1}{|F_n^d|}\sum_{\bv\in F_n^d}\pi_{S+\bv}(\mu)}_{TV}\\
    &\stackrel{(a)}{=}\dabs{\frac{1}{|F_n^d||F_{m}^d|}\sum_{\bv\in F_n^d}\sum_{\bu\in F_{m}^d} \pi_{S+\bv}(\delta_{\sigma_{n\bu}(w)})-\frac{1}{|F_n^d|}\sum_{\bv\in F_n^d}\pi_{S+\bv}(\mu)}_{TV}\\
    &\stackrel{(b)}{\leq} \frac{1}{|F_n^d|}\sum_{\bv\in F_n^d}\dabs{\frac{1}{|F_{m}^d|}\sum_{\bu\in F_{m}^d} \pi_{S+\bv}(\delta_{\sigma_{n\bu}(w)})-\pi_{S+\bv}(\mu)}_{TV}\\
    &\stackrel{(c)}{=} \frac{1}{|F_n^d|}\sum_{\bv\in F_n^d}\dabs{ \pi_{S+\bv}\parenv{\frac{1}{|F_{m}^d|}\sum_{\bu\in F_{m}^d}\delta_{\sigma_{n\bu}(w)}}-\pi_{S+\bv}(\mu)}_{TV}\\
    &= \frac{1}{|F_n^d|}\sum_{\bv\in F_n^d}\dabs{ \pi_{S+\bv}(\hfr^{F_n^d}_w)-\pi_{S+\bv}(\mu)}_{TV}\\
    &\stackrel{(d)}{\leq}\frac{1}{|F_n^d|}\sum_{\bv\in F_n^d}\dabs{ \hfr^{F_n^d}_w-\mu}_{TV} \\
    &=\dabs{ \hfr^{F_n^d}_w-\mu}_{TV}
  \end{align*}
  where:
  \begin{itemize}
  \item $(a)$ follows since
    $\pi_S(\delta_{\shift{v}(w)})=\pi_{S+\bv}(\delta_{w})$.
  \item $(b)$ follows by the triangle inequality.
  \item $(c)$ follows since $\pi_S$ is an affine map.
  \item $(d)$ follows by Lemma \ref{lem:contr}.
  \end{itemize}
  Thus, for $\epsilon>\delta$, if $\|\hfr^{F_n^d}_w-\mu\|_{TV}<\epsilon-\delta$ then $\|\fr^{S}_{w}-\hpi_S(\mu^m)\|_{TV}< \epsilon-\delta$. Therefore,
  \[
    \mathset{ w \in \Sigma^{F_{nm}^d}:~ \dabs{\fr^{S}_w- \hpi_S(\mu)}_{TV} \geq\epsilon -\delta} \subseteq
    \mathset{ w \in \Sigma^{F_{nm}^d}:~ \dabs{ \hfr^{F_n^d}_{w} - \mu}_{TV} \geq \epsilon-\delta}.
    \]
  Using the fact that $\hpi_S(\mu) \in \Bd(\Gamma)$, it follows that
  \begin{equation}
    \label{eq:fr_union}
    \mathset{ w \in \Sigma^{F_{nm}^d} ~:~ \fr^S_w \notin\sint\parenv{\Be(\Gamma)}} \subseteq
    \mathset{ w\in\Sigma^{F_{nm}^d} ~:~ \dabs{\hfr^{F_n^d}_{w}-\mu}_{TV} \geq \epsilon-\delta },
  \end{equation}
  where $\sint(\cdot)$ denotes the interior of a set, i.e.,  $\sint(\Be(\Gamma))=\mathset{\nu\in\cP(\Sigma^S) ~:~ \inf_{\mu\in\Gamma} \|\nu-\mu\|_{TV}<\epsilon}$. 
  
  If $w\in\Sigma^{F_{nm}^d}$ was randomly drawn according to $\mu^m$,
  the non-overlapping $F_n^d$-tuples are distributed i.i.d.~according
  to $\mu$. Apply Cramer's Theorem (as in \cite[Theorem 2.2.3 remark
    c]{DemZei98}) to deduce that for $\epsilon>\delta$ and for every
  $m$,
  \[\mu^m\parenv{\mathset{w\in\Sigma^{F_{nm}^d} ~:~ \dabs{ \hfr^{F_n^d}_{w} - \mu}_{TV} \geq \epsilon-\delta}} \leq
  2\exp \parenv{-m\inf_{\nu \in \cP(\Sigma^{F_n^d}):~ \|\nu - \mu\|_{TV} \geq \epsilon- \delta}H(\nu|\mu)}.\]
  Note that the function $\nu \times \mu \mapsto H(\nu|\mu)$ is continuous and strictly positive off the diagonal. Thus, for every $\epsilon >\delta$ we have
  \[c_\mu(\epsilon) \eqdef  \inf_{\nu \in \cP(\Sigma^{F_n^d}):~ \|\nu - \mu\|_{TV} \geq \epsilon - \delta} H(\nu|\mu) >0.\]
  Hence
  \begin{equation}
    \label{eq:multiforun}
    \mu^m\parenv{\mathset{w\in\Sigma^{F_{nm}^d} ~:~ \dabs{ \hfr^{F_n^d}_{w} - \mu}_{TV} \geq \epsilon-\delta}}
    \leq 2\exp\parenv{-mc_\mu(\epsilon)}.
  \end{equation}
  
  Recall that  $\cB_{nm}(\sint(\Be(\Gamma))) = \mathset{w\in\Sigma^{F_{nm}^d} ~:~
    \fr^{S}_w \in \sint(\Be(\Gamma))}$.
  By \eqref{eq:fr_union}, we have
  \begin{align}
    \label{eq:ubound_fr}
    &\mu^m\parenv{\Sigma^{F_{nm}^d}\setminus \cB_{nm}(\sint(\Be(\Gamma)))} \\ \nonumber
    &= \mu^m\parenv{\mathset{w\in\Sigma^{F_{nm}^d} ~:~ \fr^{S}_w \notin \sint(\Be(\Gamma))}}\\ \nonumber
    &\leq \mu^m\parenv{\mathset{w\in\Sigma^{F_{nm}^d} ~:~ \dabs{ \hfr^{F_n^d}_{w} - \mu}_{TV} \geq \epsilon-\delta}}.
  \end{align}
  
  Combining \eqref{eq:multiforun} and \eqref{eq:ubound_fr} we have,
  \begin{align*}
    \xi\eqdef \mu^m\parenv{\Sigma^{F_{nm}^d}\setminus \cB_{nm}\parenv{\sint(\Be(\Gamma))}}& \leq 2\exp\parenv{-mc_{\mu}(\epsilon)}.
  \end{align*}
  
  It now follows that,
  \begin{align*}
    \frac{1}{n^d} H(\mu)&=\frac{1}{(nm)^d} H(\mu^{m}) \\
    &=-\frac{1}{(nm)^d}\sum_{w\in\Sigma^{F_{nm}^d}} \mu^{m}(w)\log_2 \mu^{m}(w) \\
    &= -\frac{1}{(nm)^d}\sum_{w\in \cB_{nm}(\sint(\Be(\Gamma)))} \mu^{m}(w)\log_2 \mu^{m}(w) \\
    &\quad \; -\frac{1}{(nm)^d}\sum_{w\notin \cB_{nm}(\sint(\Be(\Gamma)))} \mu^{m}(w)\log_2 \mu^{m}(w)\\
    &\stackrel{(a)}{\leq} (1-\xi)\cdot\frac{\log_2 \abs{\cB_{nm}(\sint(\Be(\Gamma)))}}{(nm)^d}+ \xi\cdot\frac{\log_2 \abs{\Sigma^{F_{nm}^d}\setminus \cB_{nm}(\sint(\Be(\Gamma)))}}{(nm)^d}+H_2(\xi)  \\
    &\leq \frac{\log_2 \abs{\cB_{nm}(\sint(\Be(\Gamma)))}}{(nm)^d}+ 2e^{-mc_{\mu}(\epsilon)} \frac{\log_2 \abs{\Sigma}^{(nm)^d}}{(nm)^d}+H_2(\xi) \\
    &\leq \frac{1}{(nm)^d}\log_2 \abs{\cB_{nm}(\sint(\Be(\Gamma)))}
    +2e^{-mc_{\mu}(\epsilon)}\frac{1}{(nm)^d}\log_2 \abs{\Sigma}^{(nm)^d}+H_2(\xi).
  \end{align*}
  where $(a)$ follows from standard maximization of entropy arguments,
  and where $H_2(\xi)\eqdef -\xi\log_2 \xi - (1-\xi)\log_2 (1-\xi)$ is the
  binary entropy function.  This implies
  \begin{align*}
    \frac{1}{n^d}H(\mu)&=\limsupup{m}\frac{1}{n^d} H(\mu) \\
    &\leq \limsupup{m}\frac{1}{(nm)^d}\log_2 |\cB_{nm}(\sint(\Be(\Gamma)))| \\
    &\leq \limsupup{m}\frac{1}{(nm)^d}\log_2 |\cB_{nm}(\Be(\Gamma))| \\
    &\leq \ccap\parenv{\Be(\Gamma)},
  \end{align*}
  This is true for
  every $\mu\in\ocP_n(\Bd(\Gamma))$ and hence
  \[ \sup_{\mu\in\ocP_n(\Bd(\Gamma))}\frac{1}{n^d} H(\mu)\leq \ccap\parenv{\Be(\Gamma)}.\]
  Since this holds for every $n$ we have that for every $\epsilon >\delta>0$, 
  \[\hind(\Bd(\Gamma))\leq \ccap(\Be(\Gamma)).\]
  Taking the limit as $\delta\to 0$, this implies that for every
  $\epsilon >0$,
  \[\ohind(\Gamma)\leq \ccap(\Be(\Gamma)).\]
  Finally, taking the limit as $\epsilon\to 0$, it follows that
  \[\ohind(\Gamma)\leq \tcap(\Gamma).\]
\end{IEEEproof}

We summarize our results thus far by noting that for a SCS
$\Gamma\subseteq \cP(\Sigma^k)$, since $\hind (\Gamma)\leq
\hind(\God)$, Theorem \ref{indentropy1} together with \eqref{eq:summ2}
show that
\begin{align}
  \label{eq:summ3.1}
  \hind(\Gamma)\leq\hind(\Gamma^{\otimes d})&\leq\ohind(\Gamma^{\otimes d})\leq\tcap(\Gamma^{\otimes d}),\\
  \label{eq:summ3.2}
  \hind(\Gamma)\leq\hind(\Gamma^{\otimes d})&\leq\ohind(\Gamma^{\otimes d})=\ohind(\Gamma)\leq\tcap(\Gamma).
\end{align}

\section{Upper Bound on Limiting Capacity}
\label{sec:upper}

In this section we prove that if $\Gamma\subseteq\cP(\Sigma^k)$ is a
convex one-dimensional SCS and $\God$ its $d$-axial product,
then \[\limsupup{d} \tcap(\Gamma^{\otimes d})\leq \ohind(\God).\]
The main idea is to show that for any $\epsilon>0$ we are able to find
$d$ large enough for which the independence entropy is
$\epsilon$-close to $\tcap(\Gamma^{\otimes d})$. This is the main result of \cite{MeyPav14} and the proof here is an adaptation of it.

Before going into details we introduce a different form of $d$-axial product which we call the \emph{weak} $d$-axial product. For a
one dimensional SCS, $\Gamma\subseteq\cP(\Sigma^k)$,
define
\[\oGod\eqdef\mathset{\mu\in\cP(\Sigma^{F_k^d}) ~:~ \frac{1}{d}\sum_{i\in[d]}\pi_{[k]\bunt_i}(\mu) \in\Gamma},\]
and thus
\[ \cB_n\parenv{\oGod}= \mathset{w\in F_n^d ~:~ \frac{1}{d}\sum_{i\in[d]}\fr^{[k]\bunt_i}_w\in\Gamma}.\]
For the weak $d$-axial product we define,
\[\tcP_n(\oGod)\eqdef \mathset{\mu\in\parenv{\cP(\Sigma)}^{F_n^d} ~:~ \frac{1}{d}\sum_{i\in [d]} \hpi_{[k]\bunt_i}(\mu)\in\Gamma}.\]
This last definition is a relaxed version of $\God$,
since $\ocP_n(\God)$ is the set of all independent measures for which
the average of the $k$-marginals in each direction (separately)
belongs to $\Gamma$, whereas $\tcP_n(\oGod)$ is the set of all
independent measures for which the average of $k$-marginals (over all
directions) belongs to $\Gamma$.

Correspondingly, we have,
\[\thind(\oGod)\eqdef\lim_{\epsilon \to 0^+}\limsupup{n}\sup_{\mu\in\tcP_n\parenv{\Be(\Gamma)^{\boxtimes d}}}\frac{1}{n^d}H(\mu),\]
where $H(\mu)\eqdef -\sum_{w\in\Sigma^{F_n^d}}\mu(w)\log_2 \mu(w)$ is
the entropy of $\mu$.

As will become clearer later on, it will be somewhat easier to use
$\thind(\oGod)$ than $\ohind(\God)$ in this section. First, the
following lemma shows that the relaxation leading to $\thind(\oGod)$
does not affect the independence entropy.

\begin{lemma}
  \label{lem:equ2}
  Let $k\in\N$, and let $\Gamma\subseteq\cP(\Sigma^k)$ be a convex
  one-dimensional SCS, then
  \[\ohind(\Gamma)=\ohind(\God)=\thind(\oGod).\]
\end{lemma}

\begin{IEEEproof}
  By Theorem \ref{th:equ} we already know that
  $\ohind(\God)=\ohind(\Gamma)$. Thus, we are left with proving the
  last equality. Since $\Gamma$ is convex, for every $\delta >0$,
  \[\ocP_n(\Bd(\God))\subseteq\ocP_n(\Bd(\Gamma)^{\otimes d})\subseteq\tcP_n(\Bd(\Gamma)^{\boxtimes d}).\]
  Hence,
  \[\ohind(\God)\leq \thind(\oGod).\]
	
  The other direction follows essentially by using the same method as
  in the proof of Theorem \ref{th:equ}, as we now describe. Let
  $(\bv^j_i)_{i \in [n^{d-1}]}$ be an enumeration of $F_n^{j-1}\times
  \mathset{0} \times F_n^{d-j}$, i.e.,
  \[\mathset{\bv^j_0,\dots,\bv^j_{n^{d-1}-1}} =F_n^{j-1}\times \mathset{0} \times F_n^{d-j}.\]
  Fix $\delta >0$ and $\mu\in\tcP_n(\Bd(\Gamma)^{\boxtimes d})$.  For $i \in
  [n^{d-1}]$ and $j\in [d]$, define $\mu_i^j \in \cP(\Sigma^{n})$ by
  $\mu_i^j \eqdef \pi_{[n]\bunt_j+\bv^j_i}(\mu)$.  Now let
  $\hat{\mu}\in\cP(\Sigma^{dn^d})$ be the product measure satisfying
  \[\hat{\mu}(\{a\})=\prod_{j\in [d]} \prod_{i\in [n^{d-1}]} \mu_i^j(a_{in+jn^d}\dots a_{(i+1)n+jn^d-1})\]
  for every word $a=a_0\dots a_{dn^d-1}\in\Sigma^{dn^d}$.  It is
  clear that $\hat{\mu}$ is indeed a product measure, because every
  $\mu_i^j$ is also a product measure. Now,
  \begin{align*}
    \hpi_{[k]}(\hat{\mu})
    &=\frac{1}{dn^d}\sum_{i\in [dn^d]}\pi_{i+[k]}(\hat{\mu})\\
    &=\frac{1}{dn^d}\sum_{j\in [d]}\sum_{i\in [n^{d-1}]}\sum_{\ell\in [n]}\pi_{[k]+in+\ell+jn^d}(\hat{\mu}) \\
    &=\frac{1}{dn^d}\sum_{j\in [d]}\sum_{i\in [n^{d-1}]}\parenv{\sum_{\ell\in [n-k]}\pi_{[k]+in+\ell+jn^d}(\hat{\mu})+\sum_{\ell=n-k}^{n-1} \pi_{[k]+in+\ell+jn^d}(\hat{\mu})} \\
    &=\frac{1}{dn^d}\sum_{j\in [d]}\sum_{i\in [n^{d-1}]}\parenv{\sum_{\ell\in [n-k]}\pi_{[k]+\ell}(\mu_i^j)+\sum_{\ell=n-k}^{n-1} \pi_{[k]+in+\ell+jn^d}(\hat{\mu})} \\
    &=\frac{1}{dn^d}\sum_{j\in [d]}\sum_{i\in [n^{d-1}]}\parenv{\sum_{\ell\in [n]}\pi_{[k]+\ell}(\mu_i^j)-\sum_{\ell=n-k}^{n-1}\pi_{[k]+\ell}(\mu_i^j)+\sum_{\ell=n-k}^{n-1} \pi_{[k]+in+\ell+jn^d}(\hat{\mu})} \\
    &\stackrel{(a)}{=}\frac{1}{d}\sum_{j\in [d]}\hpi_{[k]\bunt_j}(\mu)-\frac{1}{dn^d}\sum_{j\in [d]}\sum_{i\in [n^{d-1}]}\sum_{\ell=n-k}^{n-1}\parenv{\pi_{[k]+\ell}(\mu_i^j)-\pi_{[k]+in+\ell+jn^d}(\hat{\mu})}.
  \end{align*}

  Recall that from the definition of $\tcP_n(\Bd(\Gamma)^{\boxtimes d})$, we have
  \[\frac{1}{d}\sum_{j\in [d]}\hpi_{[k]\bunt_j}(\mu) \in  \Bd(\Gamma).\] 
  Since
  $\parenv{\pi_{[k]+\ell}(\mu_i^j)-\pi_{[k]+in+\ell+jn^d}(\hat{\mu})}$
  is a signed measure of total variation norm at most $2$, it follows
  that
  $\hpi_{[k]}(\hat{\mu})\in\mathbb{B}_{\frac{2k}{n}+\delta}(\Gamma)$,
  so $\hat{\mu} \in
  \ocP_{dn^d}\parenv{\mathbb{B}_{\frac{2k}{n}+\delta}(\Gamma)}$.
  Hence, for every $\epsilon> \delta >0$, and every
  $\mu\in\tcP_n(\Bd(\Gamma)^{\boxtimes d})$, we can find $n_0\in\N$ such
  that for every $n>n_0$,
  $\hat{\mu}\in\ocP_{dn^d}\parenv{\Be(\Gamma)}$, and therefore,
  \[\limsupup{n}\sup_{\mu\in\tcP_n(\Bd(\Gamma)^{\boxtimes d})}\frac{1}{n^d}H(\mu) \leq \limsupup{n} \sup_{\mu\in\ocP_{dn^d}(\Be(\Gamma))}\frac{1}{dn^d}H(\mu).\]
  Thus, we obtain $\thind(\oGod)\leq \hind(\Be(\Gamma))$ for every $\epsilon>0$, and by definition it follows that 
  \[\thind(\oGod)\leq \ohind(\Gamma).\]
\end{IEEEproof}

Given a probability space $(\cX,\cF,\bP)$, denote by
$L^2(\cX,\cF,\bP,\mathbb{C}^n)$ the Hilbert space of $\cF$-measurable
functions $f:\cX\to\C^n$ satisfying
\[\dabs{f}_{L^2}^2\eqdef\int{\ip{f,f}}\mathrm{d}\bP < \infty,\] 
where  $\ip{ \cdot,\cdot }$ is the standard inner product on $\mathbb{C}^n$. 

The following lemma is based on Dirichlet's ``pigeon hole principle'' and different versions of it are used in many de-Finetti type proofs 
(see, for example, \cite{Dia1980} \cite[Lemma 4.1]{LovSze07}). 
\begin{lemma}
  \label{lem:ezer}
  Let $(\cX,\cF,\bP)$ be a probability space and let
  $\cF_0\subseteq\cF_1\subseteq\dots\subseteq\cF_m\subseteq\cF$ be a
  sequence of sub-$\sigma$-algebras. Let $f\in
  L^2(\cX,\cF,\bP,\mathbb{C}^n)$, and denote $f_j\eqdef E[f|\cF_j]$,
  the conditional expectation of $f$ with respect to the
  sub-$\sigma$-algebra $\cF_j$. Then, there exists $t\in[m]$ such that
  \[\dabs{f_{t+1}-f_{t}}^2_{L^2}\leq\frac{1}{m}\dabs{f}_{L^2}^2\]
\end{lemma}
\begin{IEEEproof}
  For every $\ell$, let $V_\ell \eqdef
  L^2(\cX,\cF_{\ell},\bP,\mathbb{C}^n)$ denote the corresponding
  sub-space of the Hilbert space $V \eqdef
  L^2(\cX,\cF,\bP,\mathbb{C}^n)$. Then $f_{\ell}$ is an orthogonal
  projection of $f$ onto $V_\ell$. Thus, $\ip{f-f_{\ell},g}_{L^2}=0$
  for every $g\in V$. Therefore,
  \[\dabs{f_m}_{L^2}^2=\sum_{\ell\in[m]}  \dabs{f_{\ell+1}-f_{\ell}}_{L^2}^2+\dabs{f_0}_{L^2}^2.\]
  Additionally, $0\leq\|f_m\|_{L^2}^2\leq\|f\|_{L^2}^2$. The result
  follows by noticing that if $m$ non-negative real numbers sum to at
  most $\|f\|_{L^2}^2$ then the value of at least one element is at
  most $\frac{1}{m}\|f\|_{L^2}^2$.
\end{IEEEproof}

Before stating the lemmas, we need the following notation. Recall that
for $k\in\N$, we defined $[k]\eqdef\mathset{0,\dots,k-1}$. We now
define $[-k]\eqdef\mathset{-1,\dots,-k}$.

\begin{lemma}
  \label{lem:ezer6}
  For every $\epsilon >0$, and any $m\in\N$, there exists $d_0\in\N$
  such that for every $d \geq d_0$, and every $n,j \in\N$, $n\geq
  j+2$, there exists a sequence of $m+1$ random subsets
  $X_0,X_1,\ldots, X_m \subseteq F_n^d$, and random variables
  $I_{t,\bv}\in[d]$, for all $t\in[m]$, $\bv\in F_n^d$, all defined on
  an appropriate probability space $(\cX,2^\cX,\bP)$, such that all
  the following hold:
  \begin{enumerate}
  \item $\bP(X_i \subseteq X_{i+1})=1$ for all $i\in[m]$.
  \item  $\bP(|X_m| \le \epsilon |F_n^d|) \geq 1- \epsilon$.
  \item For all $\bv \in F_n^d$ and $t\in[m]$, $I_{t,\bv}$ is distributed
    uniformly on $[d]$ and is independent of $X_t$. Furthermore, for every value of $X_t$, 
    \[\bP\parenv{X_t \cup ([-(j+1)]\bunt_{I_{t,\bv}}+\bv) \subseteq X_{t+1} \given X_t } \geq 1-  \epsilon.\]
  \end{enumerate}
\end{lemma}

\begin{IEEEproof}
  Choose $0<p<1$ small enough so that $1-(1-p)^{m+1}\leq\frac{\epsilon}{2}$, and
  conveniently denote $p_i\eqdef 1-(1-p)^{i+1}$. For all $i\in[m+1]$,
  consider random subsets $A_i\subseteq F_n^d$ whose coordinates are
  chosen i.i.d.~$\mathrm{Bernoulli}(p)$, i.e., $\bP(\bv\in A_i)=p$ for
  all $\bv\in F_n^d$, independently of $F_n^d\setminus
  \mathset{\bv}$. Define $X_{-1}\eqdef\emptyset$, and for all
  $i\in[m+1]$, define
  \[ X_{i}\eqdef X_{i-i}\cup A_i.\]
  Thus, $\bP(\bv\in X_i)=p_i$ for all $\bv\in F_n^d$, independently of
  $F_n^d\setminus \mathset{\bv}$.  We contend that for large enough
  $d$, the claims hold.
	
  First, it is clear that $\bP(X_i\subseteq X_{i+1})=1$ for
  $i\in[m+1]$ by construction. Second, we have
  \[ \bP\parenv{\abs{X_m}\leq \epsilon \abs{F_n^d}} \geq
  \bP\parenv{\abs{X_m}< 2p_m\abs{F_n^d}} \geq 1-e^{-2p_m^2 n^d},\]
  where the last inequality follows from Hoeffding's inequality. Since
  the right-hand side approaches $1$ when $n\geq 2$ and $d\to\infty$,
  claim $2$ holds for large enough $d$.

  We now address claim $3$. Fix $t\in[m]$ and consider $A_{t+1}$. For
  a coordinate $\bv\in F_n^d$, denote by $D(t,\bv)$ the set
  \[D(t,\bv)\eqdef \mathset{i\in [d] ~:~ \bv+[-(j+1)]\bunt_i\subseteq A_{t+1}}.\]
  If $D(t,\bv)\neq\emptyset$ then draw $I_{t,\bv}$ uniformly from
  $D(t,\bv)$. Otherwise, draw $I_{t,\bv}$ uniformly from $[d]$. Note
  that $I_{t,\bv}$ is distributed uniformly on $[d]$ since the
  distribution of $A_{t+1}$ is invariant under coordinate permutation.
  Since the coordinates in $A_{t+1}$ are chosen independently of
  $A_t,A_{t-1},\dots,A_0$ we obtain that $I_{t,\bv}$ is independent of
  $X_t$. Finally, we have
  \[
  \bP\parenv{X_t\cup ([-(j+1)]\bunt_{I_{t,\bv}}+\bv)
  \subseteq X_{t+1} \given X_t}\geq \bP(D(t,\bv)\neq\emptyset) = 1-(1-p^{j+1})^d.
  \]
  Since the right-hand side approaches $1$ as $d\to\infty$, claim $3$ holds
  for large enough $d$.
\end{IEEEproof}

If $X$ is a random variable over some probability space, we use $\bP_X$ to denote its distribution. Let
$X_0,\dots, X_{k-1}$ be random variables over the same probability
space $(\cX,2^\cX,\bP)$. We denote by $(X_0,\dots,X_{k-1})$ the vector
distributed according to their joint probability,
$\bP_{X_0,\dots,X_{k-1}}$, and denote by $(X_0\times\dots\times
X_{k-1})$ the vector distributed according to their product
probability, i.e., $\bP_{X_0\times\dots\times
  X_{k-1}}\eqdef\prod_{i\in [k]}\bP_{X_i}$.

\begin{lemma}
  \label{lem:ezer5}
  Let $\cX$ be a finite set, and $X_0,\dots,X_{k-1}$ be $k$
  random variables defined over the same probability
  space $(\cX,2^{\cX},\bP)$. Then
  \[\dabs{ \bP_{X_0,\dots, X_{k-1}}-\bP_{X_0\times\dots\times X_{k-1}}}_{TV}\leq \sum_{i=0}^{k-2} E_{X_0,\dots,X_i}\sparenv{\dabs{\bP_{X_{i+1}|X_0,\dots, X_i} - \bP_{X_{i+1}}}_{TV}}.\]
\end{lemma}

\begin{IEEEproof}
  We prove this by induction on $k$. The case of $k=1$ is trivially
  true. In the base case of $k=2$ we have,
  \begin{align}
    \label{te1}
    \dabs{ \bP_{X_0,X_1}-\bP_{X_0\times X_1} }_{TV}&=\frac{1}{2}\sum_{x_0,x_1}\abs{ \bP_{X_0,X_1}(x_0,x_1)-\bP_{X_0}(x_0)\bP_{X_1}(x_1)} \\ \nonumber
    &= \frac{1}{2}\sum_{x_0,x_1}\abs{ \bP_{X_0}(x_0)\bP_{X_1|X_0}(x_1|x_0)-\bP_{X_0}(x_0)\bP_{X_1}(x_1)}
  \end{align}
  where the sum of $x_0$ and $x_1$ is over the support of $X_0$ and $X_1$, respectively. 
  Since $\bP_{X_0}(x_0)\geq 0$ we have 
  \begin{equation}
    \label{te2}
    \frac{1}{2}\sum_{x_0,x_1\in\cX}\abs{ \bP_{X_0}(x_0) \bP_{X_1|X_0}(x_1|x_0)- \bP_{X_0}(x_0)\bP_{X_1}(x_1)}= 
    \sum_{x_0\in\cX}\bP_{X_0}(x_0)\parenv{\frac{1}{2}\sum_{x_1\in\cX}\abs{ \bP_{X_1|X_0}(x_1|x_0)- \bP_{X_1}(x_1)}}.
  \end{equation}
  Combining \eq{te1} and \eq{te2} and using the total variation
  distance definition we obtain
  \[\dabs{ \bP_{X_0,X_1}-\bP_{X_0\times X_1} }_{TV}= E_{X_0}\sparenv{\dabs{\bP_{X_1|X_0}- \bP_{X_1}}_{TV}}.\]
	
  Now assume the statement is correct for $k-1$ random variables and we show it is correct for $k$ random variables. We write
  	\[
  	\dabs{ \bP_{X_0,\dots, X_{k-1}}-\bP_{X_0\times\dots\times X_{k-1}} }_{TV}
  	=\dabs{ \bP_{X_0,\dots, X_{k-1}}-\bP_{(X_0,\dots, X_{k-2})\times X_{k-1}}+\bP_{(X_0,\dots, X_{k-2})\times X_{k-1}}-\bP_{X_0\times\dots\times X_{k-1}} }_{TV}.
  	\]
  	By applying the triangle inequality we obtain 
  	\begin{equation}\label{te43A}
  		\dabs{ \bP_{X_0,\dots, X_{k-1}}-\bP_{X_0\times\dots\times X_{k-1}} }_{TV}
  		\leq\dabs{ \bP_{X_0,\dots, X_{k-1}}-\bP_{(X_0,\dots, X_{k-2})\times X_{k-1}}}_{TV}+\dabs{\bP_{(X_0,\dots, X_{k-2})\times X_{k-1}}-\bP_{X_0\times\dots\times X_{k-1}} }_{TV}.
  	\end{equation}Considering $Y= (X_0,\ldots,X_{k-2})$ as a tuple-valued radom variable, and applying the case $k=2$ on the pair of random variables $(Y,X_{k-1})$ we have:
  	\begin{equation}\label{te43B}
  		\dabs{\bP_{X_0,\dots, X_{k-1}}-  \bP_{(X_0,\ldots,X_{k-2})\times X_{k-1}}}_{TV} \leq
  		E_{X_0,\dots, X_{k-2}}\sparenv{\dabs{\bP_{X_{k-1|(X_0,\dots, X_{k-2})}}- \bP_{X_{k-1}}}_{TV}}
  	\end{equation}
  	It is easy to check that 
  	\begin{equation}	\label{te43C}
  		\dabs{\bP_{(X_0,\dots, X_{k-2})\times X_{k-1}}-\bP_{X_0\times\dots\times X_{k-1}} }_{TV}=
  		\dabs{\bP_{X_0,\dots, X_{k-2}}-\bP_{X_0\times\dots\times X_{k-2}} }_{TV}
  	\end{equation}
  	By the induction hypothesis we have 
  	\[\dabs{\bP_{X_0,\dots, X_{k-2}}-\bP_{X_0\times\dots\times X_{k-2}} }_{TV}\leq \sum_{i=0}^{k-3}E_{X_0,\dots, X_i}\sparenv{\dabs{\bP_{X_{i+1}|X_0,\dots,X_i}-\bP_{X_{i+1}}}_{TV}}.\] 
  	Combining this with (\ref{te43A}),  (\ref{te43B}) and  (\ref{te43C}) completes the proof.

\end{IEEEproof}

For $A\subseteq F_n^d$, let $\cF_{A} \subseteq 2^{\Sigma^{F_n^d}}$
denote the $\sigma$-algebra generated by the coordinates in $A$,
namely,
\[\cF_{A} \eqdef \mathset{ \mathset{ x \in \Sigma^{F_n^d}: \pi_A(x)  \in W} ~:~ W \subseteq  \Sigma^{A} }.\]

\begin{definition}
  \label{def:etaya}
  Let $d,k,n\in\N$, $A\subseteq F_n^d$, and let
  $y\in\Sigma^{F_n^d}$. For a one-dimensional SCS, $\Gamma\subseteq
  \cP(\Sigma^k)$, and its $d$-axial-product SCS, $\God$, we define the following: 
  \begin{align*}
  	\mu^{n,d} &\;\text{ is the uniform measure over}\; \cB_{n}(\God),\\
    \mu_{y,A}&\eqdef\mu^{n,d}\parenv{\cdot \given \cF_A}(y),\\
    \eta_{y,A}&\eqdef\prod_{\bv\in F_n^d}\pi_{\mathset{\bv}}\parenv{\mu_{y,A}}.
  \end{align*}
\end{definition}

In other words, $\mu_{y,A}$ is the uniform distribution on
$\cB_n(\God)$ given whose positions in $A$ agree with $y_A$. Moreover,
$\eta_{y,A}$ is the independent version of $\mu_{y,A}$. The following statement is a particular application of \Lref{lem:ezer5} above.

\begin{lemma}
  \label{cor:ezer5}
  For every $d,n\in\N$, $i\in [d]$, and $A\subseteq F_n^d$, we have
  \[\sum_{\bv\in F_n^d}E\sparenv{\dabs{\pi_{[k]\bunt_i+\bv}(\eta_{y,A})-\pi_{[k]\bunt_i+\bv}(\mu_{y,A})}_{TV}}
  \leq 
  \sum_{\bv\in F_n^d}\sum_{j\in [k]} E\sparenv{\dabs{\pi_{\{\bv\}} (\eta_{y,A}) -\pi_{\{\bv\}}(\mu_{y,(A\cup([-j] \bunt_i + \bv)}) }_{TV}}.\]
\end{lemma}

\begin{IEEEproof}
  First note that if $k=1$ the result is immediate since all the
  summands on the left-hand side are $0$. We now examine the case of
  $k\geq 2$. For the time being, let us fix $\bv\in F_n^d$ and
  $y\in\Sigma^{F_n^d}$. We define the random variables $X_{j}$,
  $j\in[k]$, where $X_0,\dots,X_{k-1}$ is distributed according to
  $\bP_{X_0,\dots,X_{k-1}}^{y}\eqdef
  \pi_{[k]\bunt_i+\bv}(\mu_{y,A})$. In particular, each $X_j$ is
  distributed according to $\bP_{X_j}^{y}\eqdef
  \pi_{j\bunt_i+\bv}(\mu_{y,A})=\pi_{j\bunt_i+\bv}(\eta_{y,A})$. Additionally,
  $\bP_{X_0\times \dots \times X_{k-1}}^{y}=
  \pi_{[k]\bunt_i+\bv}(\eta_{y,A})$. We use the superscript $y$ to
  emphasize that these distributions depend $y$. Also for $z \in
  \Sigma^{F_n^d}$ such that $z_{A} = y_{A}$ , the conditional
  probability $\bP^y_{X_{j+1}|X_0,\dots,X_{j}}$ evaluated at $z$ is
  equal to the measure $\pi_{(j+1) \bunt_i+\bv}(\mu_{z,A\cup([j+1]
    \bunt_i + \bv)})$. By Lemma \ref{lem:ezer5}, we have
  \begin{equation}
    \label{eq321}
    \dabs{\Bbb{P}^y_{X_0,\dots,X_{k-1}}-\Bbb{P}^y_{X_0\times\dots\times X_{k-1}}}_{TV}
    \leq \sum_{j=0}^{k-2} E\sparenv{\dabs{\Bbb{P}^y_{X_{j+1}|X_0,\dots,X_{j}}- \Bbb{P}^y_{X_{i+1}}}_{TV}}.
  \end{equation}
  The expectations in the right-hand side are with respect to the
  conditioning on the random variables $X_0,\dots,X_j$.  We can
  rewrite the above equation as follows:
  \begin{equation}
    \dabs{\pi_{[k]\bunt_i+\bv}(\mu_{y,A}) - \pi_{[k]\bunt_i+\bv}(\eta_{y,A})}_{TV} 
    \leq \sum_{j=0}^{k-2} \int
    \dabs{
      \pi_{(j+1) \bunt_i+\bv}(\mu_{z,A\cup([j+1] \bunt_i + \bv)})- \pi_{(j+1)\bunt_i+\bv} (\eta_{z,A})
    }_{TV}
    \mathrm{d}\mu_{y,A}(z)
  \end{equation}
  Integrating the above inequality over $y$ with respect to
  $\mu^{n,d}$ we have:
  \begin{align*}
    &\int \dabs{\pi_{[k]\bunt_i+\bv}(\mu_{y,A}) - \pi_{[k]\bunt_i+\bv}(\eta_{y,A})}_{TV} \mathrm{d} \mu^{n,d}(y)\\
    &\quad \leq \sum_{j=0}^{k-2} \iint
    \dabs{
      \pi_{(j+1) \bunt_i+\bv}(\mu_{z,A\cup([j+1] \bunt_i + \bv)})- \pi_{(j+1)\bunt_i+\bv} (\eta_{z,A})
    }_{TV}  \mathrm{d}\mu_{y,A}(z) \mathrm{d}\mu^{n,d}(y).
  \end{align*}
  By definition of $\mu_{y,A}$ as the conditional measure, for every
  $f:\Sigma^{F_n^d} \to \mathbb{R}$ we have
  \[ \iint f(z) d\mu_{y,A}(z) \mathrm{d}\mu^{n,d}(y) = \int f(y) \mathrm{d}\mu^{n,d}(y).\]
  Writing the integeral with repect to $\mu^{n,d}$ as
  $E\sparenv{\cdot}$, we thus have
  \[E \sparenv{\dabs{ \pi_{[k]\bunt_i+\bv}(\mu_{y,A}) - \pi_{[k]\bunt_i+\bv}(\eta_{y,A})}_{TV}}
  \leq \sum_{j=0}^{k-2} E \sparenv{
    \dabs{
      \pi_{(j+1) \bunt_i+\bv}(\mu_{y,A\cup([j+1] \bunt_i + \bv)})- \pi_{(j+1)\bunt_i+\bv} (\eta_{y,A})
    }_{TV}
  }.
  \]
  Summing over all $\bv\in F_n^d$ we obtain
  \begin{equation}
    \label{eq:reflection}
    \sum_{\bv\in F_n^d}E\sparenv{\dabs{\pi_{[k]\bunt_i+\bv}(\mu_{y,A})-\pi_{[k]\bunt_i+\bv}(\eta_{y,A})}_{TV}}
    \leq \sum_{\bv\in F_n^d}\sum_{j=0}^{k-2}E\sparenv{ \dabs{\pi_{(j+1)\bunt_i+\bv}(\mu_{y,A\cup([j+1] \bunt_i + \bv)}) -\pi_{(j+1)\bunt_i+\bv} (\eta_{y,A})}_{TV}}.
  \end{equation}
  Recall that $[-j]\eqdef\mathset{-1,\dots,-j}$, hence
  \[[j+1]\bunt_i=(j+1)\bunt_i+[-(j+1)]\bunt_i.\]
  Thus, \eq{eq:reflection} can be written as 
  \begin{align}
    \label{eq:reflection2}
    &\sum_{\bv\in F_n^d}E\sparenv{\dabs{\pi_{[k]\bunt_i+\bv}(\mu_{y,A})-\pi_{[k]\bunt_i+\bv}(\eta_{y,A})}_{TV}}\\ \nonumber
    &\quad \leq \sum_{\bv\in F_n^d}\sum_{j=0}^{k-2}E\sparenv{ \dabs{\pi_{(j+1)\bunt_i+\bv}(\mu_{y,A\cup((j+1) \bunt_i + \bv+[-(j+1)]\bunt_i)}) -\pi_{(j+1)\bunt_i+\bv} (\eta_{y,A})}_{TV}}.
  \end{align}
  Since we are summing over all $\bv\in F_n^d$, and since coordinates
  are taken modulo $n$, we may write \eq{eq:reflection2} as follows,
  \begin{equation}
    \label{eq:reflection3}
    \sum_{\bv\in F_n^d}E\sparenv{\dabs{\pi_{[k]\bunt_i+\bv}(\mu_{y,A})-\pi_{[k]\bunt_i+\bv}(\eta_{y,A})}_{TV}}
    \leq \sum_{\bv\in F_n^d}\sum_{j=0}^{k-2}E\sparenv{ \dabs{\pi_{\{\bv\}}(\mu_{y,A\cup(\bv+[-(j+1)]\bunt_i)}) - \pi_{\{\bv\}} (\eta_{y,A})}_{TV}}.
  \end{equation}
  Since the total variation distance is non-negative, \eq{eq:reflection3} implies the lemma. 
\end{IEEEproof}

The following proposition, which is used to prove the main result of
this section, considers the following scenario. Assume
$y\in\Sigma^{F_n^d}$ is randomly drawn using the measure $\mu^{n,d}$,
i.e., it is drawn uniformly at random from the set of admissible words
$\cB_n(\God)$. We then study the random variable $\eta_{y,A}$ (a
measure in itself), and ask what is the probability that it resides
within the set of measures $\tcP_n\parenv{(\Be(\Gamma))^{\boxtimes
    d}}$. For convex SCSs, we prove this probability is $\epsilon$-close to $1$, assuming $d$ is sufficiently large.

\begin{proposition}
  \label{th:m1}
  Let $k\in\N$, and let $\Gamma\subseteq\cP(\Sigma^k)$ be a convex
  SCS. For any $\epsilon>0$, there exists $d_0\in\N$, such that for
  all $d\in\N$, $d\geq d_0$, $n\in\N$, $n\geq k+2$, there exists
  $A\subseteq F_n^d$, $\abs{A}\leq \epsilon n^d$, such that for
  $y\in\Sigma^{F_n^d}$ drawn randomly using the measure $\mu^{n,d}$,
  \[\mu^{n,d}\parenv{\eta_{y,A}\in\tcP_n\parenv{(\Be(\Gamma))^{\boxtimes d}}} \ge 1- \epsilon.\]
\end{proposition}

\begin{IEEEproof}
  Recall that by Definition \ref{def:etaya}, $\eta_{y,A}$ is a product
  measure, while $\mu_{y,A}$ is not necessarily so. Additionally, we
  contend that $\hpi_{[k]\bunt_i}(\mu_{y,A}) \in \Gamma$ for all $y
  \in \cB_n(\God)$, $A \subseteq F_n^d$ and $i \in [d]$. Indeed,
  \begin{align*}
    \hpi_{[k]\bunt_i}(\mu_{y,A})&=\frac{1}{|F_n^d|}\sum_{\bv\in F_n^d} \pi_{[k]\bunt_i+\bv}(\mu_{y,A}) \\
    &= \frac{1}{|F_n^d|}\sum_{\bv\in F_n^d} \frac{1}{|\pi_A^{-1}(y_A)|}\sum_{x\in\pi_A^{-1}(y_A)}\pi_{[k]\bunt_i+\bv}(\delta_{\hat{x}}) \\
    &= \frac{1}{|\pi_A^{-1}(y_A)|}\sum_{x\in\pi_A^{-1}(y_A)} \frac{1}{|F_n^d|}\sum_{\bv\in F_n^d}\pi_{[k]\bunt_i+\bv}(\delta_{\hat{x}}) \\
    &= \frac{1}{|\pi_A^{-1}(y_A)|}\sum_{x\in\pi_A^{-1}(y_A)} \pi_{[k]\bunt_i}\parenv{\frac{1}{|F_n^d|}\sum_{\bv\in F_n^d}\delta_{\sigma_{\bv}(\hat{x})}} \\
    &= \frac{1}{|\pi_A^{-1}(y_A)|}\sum_{x\in\pi_A^{-1}(y_A)} \fr_x^{[k]\bunt_i},
  \end{align*}
  where we recall that $\pi_A^{-1}(y_A)=\mathset{x\in\cB_n(\God) ~:~
    x_A=y_A}$.  Since $\fr_x^{[k]\bunt_i}\in\Gamma$ for every
  $x\in\pi_A^{-1}(y_A)$ and since $\Gamma$ is convex the contention is
  proved. Additionally, by the convexity of $\Gamma$,
  $\hpi_{[k]\bunt_i}(\mu_{y,A}) \in \Gamma$ implies
  \[ \frac{1}{d}\sum_{i\in[d]}\hpi_{[k]\bunt_i}(\mu_{y,A}) \in \Gamma.\]
	
  Draw $y\in\Sigma^{F_n^d}$ randomly using the measure
  $\mu^{n,d}$. For any $A \subseteq F_n^d$, let us denote
  \[ D_{A,y} \eqdef \dabs{ 
    \frac{1}{d}\sum_{i \in [d]}\hpi_{[k]\bunt_i}(\eta_{y,A}) -
    \frac{1}{d}\sum_{i \in [d]}\hpi_{[k]\bunt_i}(\mu_{y,A}) }_{TV}.\]
  We will use $E_y[\cdot]$ to denote expectation with respect to the
  random variable $y$ which is randomly drawn using the measure
  $\mu^{n,d}$. Denote
  \[ D_A  \eqdef E_y [D_{A,y}].\]
  To prove the theorem, it suffices to show that for any $\epsilon
  >0$, if $d$ is large enough there exists $A \subseteq F_n^d$,
  $\abs{A} \leq \epsilon n^d$, and with probability at least
  $1-\epsilon$ (with respect to $\mu^{n,d}$) we have $D_{A,y} \le
  \epsilon$.  By a standard application of the Markov inequality, it
  is sufficient to show that (under the above conditions) $D_A \leq
  \epsilon^2$.
	
  By definition, for any $A \subseteq F_n^d$ and any $y \in
  \Sigma^{F_n^d}$ we have
  \[ D_{A,y} = \dabs{ \frac{1}{d|F_n^d|}\sum_{i \in [d]} \sum_{\bv \in F_n^d} \parenv{\pi_{[k]\bunt_i+\bv}(\eta_{y,A})-
      \pi_{[k]\bunt_i+\bv}(\mu_{y,A})} }_{TV}.\]
  Applying the triangle inequality we obtain
  \[D_{A,y} \le \frac{1}{d|F_n^d|}\sum_{i \in [d]}\sum_{\bv \in F_n^d} \dabs{\pi_{[k]\bunt_i+\bv}(\eta_{y,A})- \pi_{[k]\bunt_i+\bv}(\mu_{y,A})
  }_{TV}.\] Taking the expectation, $E_y$, on both sides and using its
  linearity we get
  \[D_A \leq  \frac{1}{d|F_n^d|}\sum_{i \in [d]}\sum_{\bv \in F_n^d} E_y\sparenv{\dabs{\pi_{[k]\bunt_i+\bv}(\eta_{y,A})- \pi_{[k]\bunt_i+\bv}(\mu_{y,A})
    }_{TV}}.\]
  By \Lref{cor:ezer5} and the linearity of the expectation we obtain 
  \[D_A \leq  \sum_{j\in [k]}\frac{1}{d|F_n^d|}\sum_{i \in [d]}\sum_{\bv \in F_n^d} E_y\sparenv{\dabs{\pi_{\{\bv\}}(\eta_{y,A})- \pi_{\{\bv\}}(\mu_{y,A\cup([-(j+1)]\bunt_i+\bv)})}_{TV}}.\]

  Consider another random variable $x\in\Sigma^{F_n^d}$, also randomly
  drawn using the measure $\mu^{n,d}$. Now define $f:\Sigma^{F_n^d} \to \mathset{0,1}^{F_n^d \times \Sigma}$ by 
  \[f(x)_{(\bv,a)} \eqdef \begin{cases}
    1 & x_{\bv}=a,\\
    0 & \text{otherwise},
  \end{cases}\]
  for all $\bv\in F_n^d$ and $a\in\Sigma$. Thus, by definition we have that
  \[\pi_{\mathset{\bv}}(\eta_{y,A})(a)=\pi_{\mathset{\bv}}(\mu_{y,A})(a)=E_x\sparenv{f(x)_{(\bv,a)} \given \cF_A}(y).\]
  Since $\Sigma$ is finite we can write the total variation distance
  as a sum, and then apply the triangle inequality, which results in
  \begin{equation} 
    \label{eq:precs}
    D_{A} \leq \frac{1}{2}\sum_{ j\in [k]} \frac{1}{d|F_n^d|}\sum_{i \in [d]}\sum_{\bv \in F_n^d}
    \sum_{ a\in \Sigma} E_y \sparenv{\abs{ E_x\sparenv{f(x)_{(\bv,a)} \given \cF_A}(y) - E_x\sparenv{f(x)_{(\bv,a)} \given \cF_{A \cup([-(j+1)]\bunt_i+\bv)}}(y)}}.
  \end{equation}
  For any $j\in [k]$, viewing the expression
  \[\sum_{\bv \in F_n^d}\sum_{ a\in \Sigma} E_y\sparenv{\abs{ 
      E_x\sparenv{f(x)_{(\bv,a)} \given \cF_A}(y) - E_x\sparenv{f(x)_{(\bv,a)} \given \cF_{A \cup([-(j+1)]\bunt_i+\bv)}}(y)}}\] 
  as an inner product of a vector in $\R^{F_n^d \times \Sigma}$ whose $(\bv,a)$'th coordinate is equal to
\[E_y\sparenv{\abs{ E_x\sparenv{f(x)_{(\bv,a)} \given \cF_A}(y) -
E_x\sparenv{f(x)_{(\bv,a)} \given \cF_{A\cup([-(j+1)] \bunt_i+\bv)}}(y)}}\] 
  and $\bone$, we may apply
  Cauchy-Schwarz (C.S) inequality and obtain
  \begin{align}
    \label{eq:cs}
    &\sum_{(\bv,a)\in F_n^d\times\Sigma} E_y\sparenv{\abs{ 	E_x\sparenv{f(x)_{(\bv,a)} \given \cF_A}(y) -E_x\sparenv{f(x)_{(\bv,a)} \given \cF_{A \cup([-(j+1)]\bunt_i+\bv)}}(y)}}\cdot 1 \\ \nonumber
    &\stackrel{\text{C.S}}{\leq} \sqrt{\parenv{\sum_{(\bv,a)\in F_n^d\times\Sigma} \parenv{E_y\sparenv{\abs{ E_x\sparenv{f(x)_{(\bv,a)} \given \cF_A}(y) - E_x\sparenv{f(x)_{(\bv,a)} \given \cF_{A \cup([-(j+1)]\bunt_i+\bv)}}(y)}}}^2} \parenv{\sum_{(\bv,a)\in F_n^d\times\Sigma} 1^2}} \\ \nonumber
    &= \sqrt{\parenv{\sum_{(\bv,a)\in F_n^d\times\Sigma} \parenv{E_y\sparenv{\abs{E_x\sparenv{f(x)_{(\bv,a)} \given \cF_A}(y) - E_x\sparenv{ f(x)_{(\bv,a)} \given \cF_{A \cup([-(j+1)]\bunt_i+\bv)}}(y)}}}^2} \cdot\abs{F_n^d}\cdot\abs{\Sigma}}.
  \end{align}
  Thus, combining \eq{eq:precs} and \eq{eq:cs} we have	
  \[D_A \leq \frac{\sqrt{\abs{\Sigma}}}{2d\sqrt{\abs{F_n^d}}} \sum_{j\in [k]}\sum_{i\in [d]} \sqrt{\sum_{\bv \in F_n^d}
    \sum_{ a\in \Sigma} \parenv{E_y\sparenv{\abs{
          E_x\sparenv{f(x)_{(\bv,a)} \given \cF_A}(y) -
          E_x\sparenv{f(x)_{(\bv,a)} \given \cF_{A
              \cup([-(j+1)]\bunt_i+\bv)}}(y)}}}^2}.\] 
	Using the fact that $(E[\abs{X}])^2\leq E[X^2]$ (again, by C.S), we have
  \begin{equation}
    \label{eq:ma1}
    D_A \leq \frac{\sqrt{\abs{\Sigma}}}{2d\sqrt{\abs{F_n^d}}} \sum_{j\in [k]}\sum_{i\in [d]} \sqrt{\sum_{\bv \in F_n^d} \sum_{ a\in \Sigma} \parenv{E_y\sparenv{\parenv{ E_x\sparenv{f(x)_{(\bv,a)} \given \cF_A}(y) - E_x\sparenv{f(x)_{(\bv,a)} \given \cF_{A \cup([-(j+1)]\bunt_i+\bv)}}(y)}^2}}}.
  \end{equation}
	
  Choose $m$ large enough such that
  $\frac{1}{\sqrt{m}}\leq\frac{\epsilon^2}{k\sqrt{\abs{\Sigma}}}$ and
  denote $\epsilon_0= \frac{\epsilon^2}{k\abs{\Sigma}}$. Now let
  $\bP,I_{t,\bv},X_0, X_1,\ldots, X_m$ be as given by Lemma \ref{lem:ezer6} with $n>k+2$ and with
  $\epsilon_0$ and obtain $d_0$. From here on, assume $d\geq d_0$. Let
  $\E$ denote the expectation with respect to $\bP$.
	
  First, from \eq{eq:ma1} we may bound $D_{X_t}$, for any $t\in[m+1]$, by 
  \begin{equation}
    \label{eq:ma1a}
    D_{X_t} \leq \frac{\sqrt{\abs{\Sigma}}}{2\sqrt{\abs{F_n^d}}}\sum_{ j\in [k]} \frac{1}{d}\sum_{i \in [d]}\sqrt{\sum_{\bv \in F_n^d}\sum_{ a\in \Sigma} E_y\sparenv{\parenv{ E_x\sparenv{ f(x)_{(\bv,a)} \given \cF_{X_t}}(y) - E_x\sparenv { f(x)_{(\bv,a)} \given \cF_{X_{t}\cup ([-(j+1)]\bunt_i+\bv)}}(y)}^2}}.
  \end{equation}
  By the properties of $X_t$ and $X_{t+1}$ given in Lemma
  \ref{lem:ezer6}, for every $\bv \in F_n^d$ there is a random
  variable $I_{t,\bv}$ independent of $X_t$ and distributed
  uniformly on $[d]$ so that $\bP(X_t \cup ([-(j+1)]\bunt_{I_{t,\bv}} +\bv)
  \subseteq X_{t+1} \mid X_t) \ge 1- \epsilon_0$. Denote
  \[X_{t,\bv} \eqdef X_t \cup ([-(j+1)]\bunt_{I_{t,\bv}} +\bv).\]
  Since $I_{t,\bv}$ is independent of $X_t$ we have
  \begin{align}
    \label{eq:ma2}
    &\E\sparenv{\sqrt{\sum_{\bv \in F_n^d}\sum_{ a\in \Sigma} E_y\sparenv{\parenv{ E_x\sparenv{ f(x)_{(\bv,a)} \given \cF_{X_t}}(y)-E_x\sparenv { f(x)_{(\bv,a)} \given \cF_{X_{t,\bv}}}(y)}^2}} \given X_t} \\ \nonumber
    &=	\frac{1}{d}\sum_{i \in [d]}  \sqrt{\sum_{\bv \in F_n^d}\sum_{ a\in \Sigma} E_y\sparenv{\parenv{ E_x\sparenv{ f(x)_{(\bv,a)} \given \cF_{X_t}}(y) - E_x\sparenv { f(x)_{(\bv,a)} \given \cF_{X_{t}\cup ([-(j+1)]\bunt_i +\bv)}}(y)}^2}}. 
  \end{align}
  From \eq{eq:ma1a} and \eq{eq:ma2} we obtain 
  \begin{align}
    \label{eq:ma2a}
    D_{X_t} &\leq \frac{\sqrt{\abs{\Sigma}}}{2\sqrt{\abs{F_n^d}}}\sum_{ j\in [k]} \E\sparenv{\sqrt{\sum_{\bv \in F_n^d}\sum_{ a\in \Sigma} E_y\sparenv{\parenv{ E_x\sparenv{ f(x)_{(\bv,a)} \given \cF_{X_t}}(y) - E_x\sparenv { f(x)_{(\bv,a)} \given \cF_{X_{t,\bv}}}(y)}^2}} \given X_t}.
  \end{align}
		
  Since we may view $E_x\sparenv{f(x)_{(\bv,a)} \given \cF_{X_t}}$ as
  the orthogonal projection of $f(x)_{(\bv,a)}$ on
  $L^2(F_n^d,\cF_{X_t},\mu^{n,d},\R)$, if $X_{t,\bv}\subseteq X_{t+1}$
  we have
  \begin{align*}
    &E_y\sparenv{\parenv{E_x\sparenv{f(x)_{(\bv,a)} \given \cF_{X_t}}(y) - E_x\sparenv{f(x)_{(\bv,a)} \given \cF_{X_{t,\bv}}}(y) }^2} \\
    &\leq E_y\sparenv{\parenv{E_x\sparenv{f(x)_{(\bv,a)} \given \cF_{X_t}}(y) - E_x\sparenv{f(x)_{(\bv,a)} \given \cF_{X_{t+1}}}(y) }^2}.
  \end{align*}
  Otherwise, if $X_{t,\bv}\nsubseteq X_{t+1}$, then 
  \[E_y\sparenv{\parenv{E_x\sparenv{f(x)_{(\bv,a)} \given \cF_{X_t}}(y) - E_x\sparenv{f(x)_{(\bv,a)} \given \cF_{X_{t,\bv}}}(y) }^2}\leq 1.\] 
  By the properties of $X_t$ given in \Lref{lem:ezer6}, we have that  
  \[\bP(X_{t,\bv}\subseteq X_{t+1}\mid X_t)\geq 1-\epsilon_0.\]
  Thus, 
  \begin{align}
    \label{eq:ma3}
    &\E\sparenv{\sqrt{\sum_{\bv \in F_n^d}\sum_{ a\in \Sigma} E_y\sparenv{\parenv{ E_x\sparenv{ f(x)_{(\bv,a)} \given \cF_{X_t}}(y) - E_x\sparenv { f(x)_{(\bv,a)} \given \cF_{X_{t,\bv}}}(y)}^2}}\given X_t} \\ \nonumber
    &\leq (1-\epsilon_0)\E\sparenv{\sqrt{\sum_{\bv \in F_n^d}\sum_{ a\in \Sigma} E_y\sparenv{\parenv{ E_x\sparenv{ f(x)_{(\bv,a)} \given \cF_{X_t}}(y) - E_x\sparenv { f(x)_{(\bv,a)} \given \cF_{X_{t+1}}}(y)}^2}}\given X_t}+\epsilon_0\sqrt{\abs{\Sigma}\abs{F_n^d}} \\ \nonumber
    &\leq \E\sparenv{\sqrt{\sum_{\bv \in F_n^d}\sum_{ a\in \Sigma} E_y\sparenv{\parenv{ E_x\sparenv{ f(x)_{(\bv,a)} \given \cF_{X_t}}(y) - E_x\sparenv { f(x)_{(\bv,a)} \given \cF_{X_{t+1}}}(y)}^2}}\given X_t}+\epsilon_0\sqrt{\abs{\Sigma}\abs{F_n^d}}.
  \end{align}
  From \eq{eq:ma2a} and \eq{eq:ma3} we obtain 
  \begin{align}
    \label{eq:ma3a}
    D_{X_t} &\leq \frac{\sqrt{\abs{\Sigma}}}{2\sqrt{\abs{F_n^d}}}\sum_{ j\in [k]} \E\sparenv{\sqrt{\sum_{\bv \in F_n^d}\sum_{ a\in \Sigma} E_y\sparenv{\parenv{ E_x\sparenv{ f(x)_{(\bv,a)} \given \cF_{X_t}}(y) - E_x\sparenv { f(x)_{(\bv,a)} \given \cF_{X_{t+1}}}(y)}^2}}\given X_t} \\ \nonumber
    &\quad\  + \frac{\sqrt{\abs{\Sigma}}}{2\sqrt{\abs{F_n^d}}}\sum_{ j\in [k]}\epsilon_0\sqrt{\abs{\Sigma}\abs{F_n^d}} \\ \nonumber
    &= \frac{k\sqrt{\abs{\Sigma}}}{2\sqrt{\abs{F_n^d}}}\E\sparenv{\sqrt{\sum_{\bv \in F_n^d}\sum_{ a\in \Sigma} E_y\sparenv{\parenv{ E_x\sparenv{ f(x)_{(\bv,a)} \given \cF_{X_t}}(y) - E_x\sparenv { f(x)_{(\bv,a)} \given \cF_{X_{t+1}}}(y)}^2}}\given X_t} + \frac{k\abs{\Sigma}\epsilon_0}{2}.
  \end{align}

  Observe that viewing $f$ as a random variable with respect to
  $\mu^{n,d}$ we have $\dabs{f}_2 = \sqrt{\abs{F_n^d}}$. Note also
  that
  \begin{equation}
    \label{eq:ma4}
    \dabs{ E\sparenv{f \given \cF_{X_t}} - E\sparenv{f \given \cF_{X_{t+1}}}}_2=\sqrt{\sum_{\bv \in F_n^d}\sum_{ a\in \Sigma}	E_y\sparenv{\parenv{E_x\sparenv{f(x)_{(\bv,a)} \given \cF_{X_t}}(y) - E_x\sparenv{f(x)_{(\bv,a)} \given \cF_{X_{t+1}}}(y) }^2}}.
  \end{equation}
  From \eq{eq:ma3a} and \eq{eq:ma4} we obtain that for every $t\in [m]$, 
  \begin{align}
    \label{eq:fin}
    D_{X_t}\leq \frac{k\sqrt{\abs{\Sigma}}}{2\sqrt{\abs{F_n^d}}}\E\sparenv{\dabs{ E[f \given \cF_{X_t}] - E[f \given \cF_{X_{t+1}}]}_2\given X_t}+\frac{k\abs{\Sigma}\epsilon_0}{2}.
  \end{align}
  Note that the probability that a random variable is greater or equal to its expectation is always strictly positive.  Because $X_t$ takes only finitely many values, this means that for every
  $t\in [m]$, for every realization of $X_t$, denoted as $\chi_t$, there
  exists a realization of $X_{t+1}$, denoted as $\chi_{t+1}=\chi_{t+1}(\chi_t)$
  such that $\bP\parenv{X_{t+1}=\chi_{t+1}\given X_t}>0$ and
  \[\E\sparenv{\dabs{E\sparenv{f\given \cF_{X_t}}-E\sparenv{f\given \cF_{X_{t+1}}}}_2\given X_t}\leq\dabs{E\sparenv{f\given \cF_{\chi_t}}-E\sparenv{f\given \cF_{\chi_{t+1}}}}_2.\]
  Together with \eq{eq:fin} we obtain
  \begin{align}
    \label{eq:fin1}
    D_{\chi_t}\leq \frac{k\sqrt{\abs{\Sigma}}}{2\sqrt{\abs{F_n^d}}}\dabs{ E\sparenv{f \given \cF_{\chi_t}} - E\sparenv{f \given\cF_{\chi_{t+1}}}}_2+\frac{k\abs{\Sigma}\epsilon_0}{2}.
  \end{align}
  Since \eq{eq:fin1} holds for every $t$, we obtain that there exists
  a sequence $(\chi_t)_{t\in [m+1]}$ of realizations of $(X_t)_{t\in
    [m+1]}$ with positive probabilities, such that for every $t\in
  [m]$,
  \begin{align}
    \label{eq:fin2}
    D_{\chi_t}\leq \frac{k\sqrt{\abs{\Sigma}}}{2\sqrt{\abs{F_n^d}}}\dabs{ E\sparenv{f \given \cF_{\chi_t}} - E\sparenv{f \given \cF_{\chi_{t+1}}}}_2+\frac{k\abs{\Sigma}\epsilon_0}{2}.
  \end{align}
  From \Lref{lem:ezer}, there exists $t\in [m]$ such that 
  \begin{equation}
    \label{eq:ma4a}
    \dabs{ E\sparenv{f \given \cF_{\chi_t}} - E\sparenv{f \given \cF_{\chi_{t+1}}} }_2^2\leq \frac{1}{m}\dabs{f}_2^2.
  \end{equation}
  Combining \eq{eq:ma4a} with \eq{eq:fin2} we obtain that there exists
  $t\in [m]$ such that
  \begin{align*}
    D_{\chi_t}&\leq \frac{k\sqrt{\abs{\Sigma}}}{2\sqrt{\abs{F_n^d}}}\frac{1}{\sqrt{m}}{\dabs{f}_2}+\frac{k\abs{\Sigma}\epsilon_0}{2} \\
    &= \frac{k\sqrt{\abs{\Sigma}}}{2\sqrt{m}}+ \frac{k\abs{\Sigma}\epsilon_0}{2}.
  \end{align*}
  Taking $A=\chi_t$, and recalling our choice of
  $\epsilon_0=\frac{\epsilon^2}{k\abs{\Sigma}}$ and
  $\frac{1}{\sqrt{m}}\leq \frac{\epsilon^2}{k\sqrt{\abs{\Sigma}}}$, we
  obtain that
  \begin{align*}
    D_{A}&\leq \frac{k\sqrt{\abs{\Sigma}}}{2\sqrt{m}}+ \frac{k\abs{\Sigma}\epsilon_0}{2} \leq \epsilon^2,
  \end{align*}
  which completes the proof.
\end{IEEEproof}

We have reached the main result of this section. We show that capacity of a convex $d$-axial product is arbitrarily close to the independence entropy, as the dimension grows.

\begin{theorem}
  \label{th:main2}
  Let $k\in\N$, and let $\Gamma \subseteq \cP(\Sigma^k)$ be a convex
  one-dimensional SCS. Then
  \[\limsupup{d}\tcap(\God) = \ohind(\Gamma).\]
\end{theorem}

\begin{IEEEproof}
  First note that $\limsupup{d}\tcap(\God) \geq \ohind(\Gamma)$ by applying Theorem \ref{indentropy1} to $\God$ for every $d$ and taking $d\to \infty$ on both sides. 
  For the other direction, fix $\epsilon_0>0$ and choose
  \[0<\epsilon< \min\mathset{\frac{\epsilon_0}{2\log_2 \abs{\Sigma}},1}, \qquad 0<\delta <\frac{\epsilon}{2}.\]
  Replace $\Gamma$ by $\Bd(\Gamma)$ in Definition \ref{def:etaya} and
  denote the resulting measures by $\mu^{n,d}_\delta$,
  $\mu_{y,A}^{\delta}$, and $\eta_{y,A}^{\delta}$.

  Recall that for a measure $\mu$ and a $\sigma$-algebra $\cF$,
  \begin{equation}
    \label{eq:def1}
    H\parenv{\mu \given \cF}\eqdef E \sparenv{H(\mu\parenv{\cdot\given \cF})}=
    \int H\parenv{\mu\parenv{\cdot\given \cF}}(x) \mathrm{d}\mu(x).
  \end{equation}
  In other words, $H(\mu \mid \cF)$ is
  the expected entropy of the conditional measure
  $\mu\parenv{\cdot\given \cF}$. Also recall that for $A\subseteq
  F_n^d$, $\pi_A(\mu^{n,d}_\delta)$ denotes the $A$-marginal of
  $\mu^{n,d}_\delta$, and that $\cF_A$ denotes the $\sigma$-algebra
  generated by the coordinates in $A$.  We have that
  \[H(\mu^{n,d}_\delta)=H(\pi_A(\mu^{n,d}_\delta))+H\parenv{\mu^{n,d}_\delta \given \cF_A}.\]

  By \Pref{th:m1}, for any $n\in\N$, $n\geq k+2$, there exists $d_0
  \in \N$, such that for every $d \geq d_0$, there exists $A \subseteq
  \Sigma^{F_n^d}$, $\abs{A}\leq \epsilon n^d$, such that,
  \[\mu^{n,d}_\delta \parenv{\eta_{y,A}^\delta\in\tcP_n\parenv{(\Be(\Bd(\Gamma)))^{\boxtimes d}}} \geq 1- \epsilon > 0.\]
  In particular, there exists a word $y\in\Sigma^{F_n^d}$ such that
  $\eta_{y,A}\in\tcP_n\parenv{(\Be (\Bd(\Gamma)))^{\boxtimes
      d}}$. Since clearly
  \[H(\pi_A(\mu^{n,d}_\delta)) \leq \log_2 \abs{\Sigma^A},\]
  by combining the above we have
  \begin{equation}
    \label{eq:cond_entropy_bound}
    H(\mu^{n,d}_\delta) \leq H\parenv{\mu^{n,d}_\delta\given \cF_A}+ \epsilon n^d \log_2 \abs{\Sigma}.
  \end{equation}
  Because the joint entropy of a finite set of random variables is
  bounded from above by the sum of their entropies (and the same
  statement holds for conditional entropy), we have:
  \[ H\parenv{\mu^{n,d}_\delta \given \cF_A} \leq \sum_{ \bv \in F_n^d}  H\parenv{\pi_{\{\bv\}}(\mu^{n,d}_\delta)\given \cF_A}.\]
  By definition of the random measure $\eta_{y,A}^\delta$ and from
  \eq{eq:def1}, we have
  \[ H\parenv{\pi_{\mathset{\bv}}(\mu^{n,d}_\delta)\given \cF_A}= \sum_{y\in\Sigma^{F_n^d}} H\parenv{\pi_{\mathset{\bv}}(\eta_{y,A}^\delta)} \mu^{n,d}_\delta(y).\]
  Thus, 
  \[H\parenv{\mu^{n,d}_\delta\given \cF_A} \leq \sum_{ \bv \in F_n^d} \sum_{y\in\Sigma^{F_n^d}} H\parenv{\pi_{\{\bv\}}(\eta^\delta_{y,A})} \mu^{n,d}_\delta(y).\]
  Now, since $\eta_{y,A}$ is a product measure, we have
  \[H(\eta_{y,A}^\delta) = \sum_{\bv \in F_n^d}H\parenv{\pi_{\{\bv\}} (\eta_{y,A}^\delta)}.\]
  It follows that, 
  \begin{equation}
    \label{eq:na1}
    H\parenv{\mu^{n,d}_\delta\given \cF_A} \leq \sum_{y\in\Sigma^{F_n^d}} H\parenv{\eta_{y,A}^\delta}\mu^{n,d}_\delta(y).
  \end{equation}
	
  Let us conveniently use $p$ to denote the value
  \[p \eqdef  \mu^{n,d}_\delta\parenv{\eta_{y,A}^\delta \in  \tcP_n\parenv{(\Be(\Bd(\Gamma)))^{\boxtimes d}}},\] 
  and recall that $p\geq 1-\epsilon>0$. Then
  \begin{align}
    \label{eq:na2}
    \sum_{y\in \Sigma^{F_n^d}} H(\eta_{y,A}^\delta)\mu^{n,d}_\delta(y) \leq p\cdot \sup_{\eta\in \tcP_n \parenv{(\Be(\Bd(\Gamma)))^{\boxtimes d}}} H(\eta)+(1-p)\cdot\log_2\abs{\Sigma^{F_n^d}}.
  \end{align}
  Using the fact that $p \geq 1- \epsilon>0$ combined with \eq{eq:na1}
  and \eq{eq:na2}, it follows that
  \begin{align} 
    \label{eq:na3}
    H\parenv{\mu^{n,d}_\delta \given \cF_A} &\leq p\cdot \sup_{\eta\in \tcP_n\parenv{(\Be(\Bd(\Gamma)))^{\boxtimes d}}}H(\eta) + (1-p)n^d\cdot \log_2 \abs{\Sigma} \\ \nonumber
    &\leq \sup_{\eta\in \tcP_n\parenv{(\Be(\Bd(\Gamma)))^{\boxtimes d}}}H(\eta) + \epsilon n^d \log_2 \abs{\Sigma}.
  \end{align}
  Combining \eq{eq:na3} with \eq{eq:cond_entropy_bound} we obtain
  \begin{equation*}
    \frac{1}{n^d}H\parenv{\mu^{n,d}_\delta} \leq  \frac{1}{n^d}\sup_{\eta\in \tcP_n\parenv{(\Be(\Bd(\Gamma)))^{\boxtimes d}}}H(\eta) + 2\epsilon\log_2 \abs{\Sigma}.
  \end{equation*}
  By our choice of $\epsilon$, we have $\epsilon+\delta \leq
  \epsilon_0$, hence $(\Be(\Bd(\Gamma)))^{\otimes d}\subseteq
  (\Bbb{B}_{\epsilon_0}(\Gamma))^{\otimes d}$, as well as
  \begin{equation}\label{eq:H_mu_n_d_ent_bound_2}
    \frac{1}{n^d}H\parenv{\mu^{n,d}_\delta} \leq  \frac{1}{n^d}\sup_{\eta\in \tcP_n\parenv{(\Bbb{B}_{\epsilon_0}(\Gamma))^{\boxtimes d}}}H(\eta) + \epsilon_0.
  \end{equation}
  Since $\mu^{n,d}_\delta$ is the uniform measure on
  $\cB_n(\Bd(\Gamma)^{\otimes d})$,
  \[H(\mu^{n,d}_\delta)=\log_2 \abs{\cB_n(\Bd(\Gamma)^{\otimes d})}.\] 
  Thus, 
  \[\frac{1}{n^d}\log_2 \abs{\cB_n\parenv{(\Bd(\Gamma))^{\otimes d}}} \leq  \frac{1}{n^d}\sup_{\eta\in \tcP_n\parenv{(\Bbb{B}_{\epsilon_0}(\Gamma))^{\boxtimes d}}}H(\eta) + \epsilon_0.\]
  Taking $\limsup_{n \to \infty}$ we obtain
  \[
  \ccap\parenv{(\Bd(\Gamma))^{\otimes d}} \leq \limsupup{n}\frac{1}{n^d}\sup_{\eta\in \tcP_n\parenv{(\Bbb{B}_{\epsilon_0}(\Gamma))^{\boxtimes d}}}H(\eta) + \epsilon_0 .
  \]
  Since
  \[  \Bd \parenv{\God } \subseteq \parenv{\Bd(\Gamma)}^{\otimes d},\]
  we have
  \[
  \ccap\parenv{\Bd\parenv{(\Gamma)^{\otimes d}}} \leq \ccap\parenv{(\Bd(\Gamma))^{\otimes d}}\leq \limsupup{n}\frac{1}{n^d}\sup_{\eta\in \tcP_n\parenv{(\Bbb{B}_{\epsilon_0}(\Gamma))^{\boxtimes d}}}H(\eta) + \epsilon_0 .
  \]
  Taking $\lim_{\delta\to 0^+}$, we get
  \begin{equation}
  \label{eq351}
  \tcap\parenv{\God} \leq \limsupup{n}\frac{1}{n^d}\sup_{\eta\in \tcP_n\parenv{(\Bbb{B}_{\epsilon_0}(\Gamma))^{\boxtimes d}}}H(\eta) + \epsilon_0 .
  \end{equation}
  At this point we take a slight detour. For $\xi>0$,
  $\Bbb{B}_{\epsilon_0}(\Gamma)\subseteq\Bbb{B}_{\xi}\parenv{\Bbb{B}_{\epsilon_0}(\Gamma)}$
  and hence we have
  \begin{align*}
    \limsupup{n}\frac{1}{n^d}\sup_{\eta\in \tcP_n\parenv{(\Bbb{B}_{\epsilon_0}(\Gamma))^{\boxtimes d}}}H(\eta) + \epsilon_0 &\leq \limsup_{\xi\to 0^+} \limsupup{n}\frac{1}{n^d}\sup_{\eta\in \tcP_n\parenv{(\Bbb{B}_{\xi}(\Bbb{B}_{\epsilon_0}(\Gamma)))^{\boxtimes d}}}H(\eta) + \epsilon_0\\
    &\overset{(a)}{=}\thind\parenv{\parenv{\Bbb{B}_{\epsilon_0}(\Gamma)}^{\boxtimes d}}+\epsilon_0\\
    &\overset{(b)}{=} \ohind\parenv{\parenv{\Bbb{B}_{\epsilon_0}(\Gamma)}}+\epsilon_0,
  \end{align*}
  where $(a)$ follows by definition, and $(b)$ follows by Lemma \ref{lem:equ2}.
  Substituting this in \eq{eq351} and taking $d\to \infty$ we obtain 
  \begin{equation}
    \label{eq352}
    \limsupup{d}\tcap\parenv{\God} \leq \ohind\parenv{\parenv{\Bbb{B}_{\epsilon_0}(\Gamma)}}+\epsilon_0.
  \end{equation}
  Note that since $\Gamma$ is convex we have that for $\epsilon_1>0$,
  $\Bbb{B}_{\epsilon_1}\parenv{\Bbb{B}_{\epsilon_0}(\Gamma)}=\Bbb{B}_{\epsilon_1+\epsilon_0}(\Gamma)$. Therefore,
  by the definition of limit we have
\[\limsup_{\epsilon_0\to 0^+}\limsup_{\epsilon_1\to 0^+} \limsupup{n}\frac{1}{n}\sup_{\eta\in \ocP_n\parenv{(\Bbb{B}_{\epsilon_1}(\Bbb{B}_{\epsilon_0}(\Gamma)))}}H(\eta)= \limsup_{\epsilon_0\to 0^+} \limsupup{n}\frac{1}{n}\sup_{\eta\in \ocP_n\parenv{(\Bbb{B}_{\epsilon_0}(\Gamma))}}H(\eta).\]
Therefore, taking the limit as $\epsilon_0\to 0$ in \eq{eq352} we
obtain
\[ \limsupup{d}\tcap\parenv{\God} \leq \ohind(\Gamma).\]
\end{IEEEproof}

\section{Discussion}
\label{sec:conc}

Our initial motivation behind this work is to approximate the capacity of multidimensional SCSs using ``meaningful'' expressions. The main challenges were defining exactly what is the capacity of multidimensional SCSs, and obtaining the connections between the capacity and the independence entropy. Our approach, which uses the independence entropy, extends previous combinatorial works
\cite{PooChaMar06,LouMarPav13,MeyPav14}, which apply only to fully
constrained systems. At the core of our results, for
$\Gamma\subseteq\cP(\Sigma^k)$ and its axial product $\God$, by
Theorem \ref{th:equ} and Theorem \ref{indentropy1} that
\[ \ohind(\Gamma)\leq \tcap(\God).\]
Thus, the problem of bounding the capacity of a $d$-dimensional
axial-product SCS is simplified by having to consider only product
measures, which are much easier to handle. Moreover, any number of
dimensions $d$, may be reduced via this bound to the one-dimensional
case. This bound is asymptotically tight, as together with Theorem
\ref{th:main2}, for convex $\Gamma$,
\[\limsupup{d}\tcap(\God) = \ohind(\Gamma).\]
It also appears that the capacity $\tcap$, and
independence entropy $\ohind$, are robust generalizations of their
one-dimensional combinatorial counterparts.

\begin{figure}[ht]
  \begin{center}
    
    \begin{tikzpicture}[scale=1]
      
      \draw[fill] (0,0) circle (.05cm) node[right] {\ $\ohind(\Gamma)=\ohind(\God)=\thind(\oGod)$};
      \draw[fill] (-2,2) circle (.05cm) node[left] {$\tcap(\Gamma)$\ \ \ };
      \draw[fill] (-4,1) circle (.05cm) node[left] {$\ccap(\Gamma)$};
      \draw[fill] (2,2) circle (.05cm) node[right] {\ \ $\tcap(\God)$};
      \draw[fill] (4,1) circle (.05cm) node[right] {$\ccap(\God)$};
      \draw[fill] (0,-2) circle (.05cm) node[right] {$\hind(\God)$};
      \draw[fill] (0,-4) circle (.05cm) node[right] {$\hind(\Gamma)$};
      
      \draw (0,0) -- (-2,2);
      \draw (-2,2) -- (-4,1);
      \draw (0,0) -- (2,2);
      \draw (2,2) -- (4,1);
      \draw (0,0) -- (0,-2);
      \draw (0,-2) -- (0,-4);

      \node[text width=1cm, anchor=west, right] at (-3.25,1.25) {(a)};
      \node[text width=1cm, anchor=west, right] at (-1.1,1.25) {(b)};
      \node[text width=1cm, anchor=west, right] at (0.5,1.25) {(c)};
      \node[text width=1cm, anchor=west, right] at (2.5,1.25) {(d)};
      \node[text width=1cm, anchor=west, right] at (-0.75,-1) {(e)};
      \node[text width=1cm, anchor=west, right] at (-0.75,-3) {(f)};
      
    \end{tikzpicture}
    
  \end{center}
  \caption{The Hasse diagram for a convex one-dimensional $\Gamma\subseteq\cP(\Sigma^k)$ and its $d$-axial-product SCS $\God$, where (a) and (d) follow from \eqref{eq:ccaptcap}, (b) and (c) follow from Theorem \ref{indentropy1}, (e) follows from \eqref{eq:hindohind}, and (f) follows from Lemma \ref{lem:ind_axial}. }
  \label{fig:axial}
\end{figure}
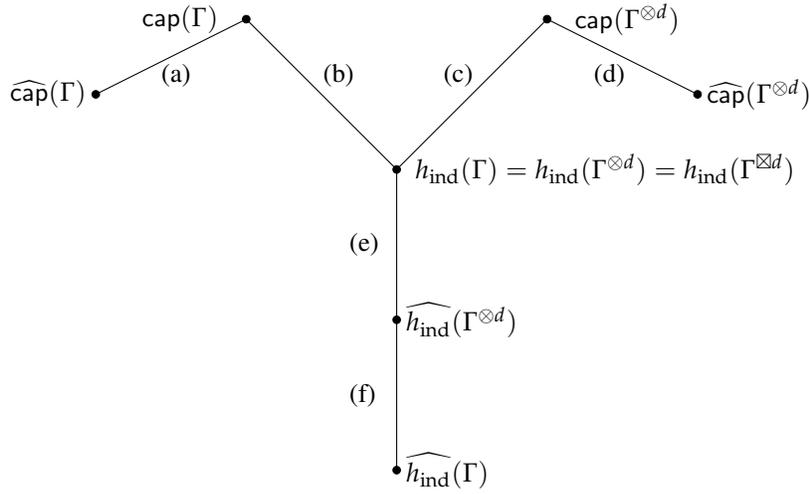

The paper contains many connections between the various capacities and
entropies. Figure \ref{fig:general} shows the Hasse diagram for the
bounds pertaining to general $d$-dimensional
$\Gamma\subseteq\cP(\Sigma^{F_k^d})$. In the case of a convex
one-dimensional $\Gamma\subseteq\cP(\Sigma^k)$ and its
$d$-axial-product SCS $\God$, a more elaborate Hasse diagram emerges,
which is shown in Figure \ref{fig:axial}.

\begin{figure}[ht]
  \begin{center}
    
    \begin{tikzpicture}[scale=1]
      
      \draw[fill] (0,0) circle (.05cm) node[right] {\ $\tcap(\Gamma)$};
      \draw[fill] (-2,-2) circle (.05cm) node[left] {$\ccap(\Gamma)$};
      \draw[fill] (2,-2) circle (.05cm) node[right] {$\ohind(\Gamma)$};
      \draw[fill] (2,-4) circle (.05cm) node[right] {$\hind(\Gamma)$};
      
      \draw (0,0) -- (-2,-2);
      \draw (0,0) -- (2,-2);
      \draw (2,-2) -- (2,-4);

      \node[text width=1cm, anchor=west, right] at (-1.75,-0.75) {(a)};
      \node[text width=1cm, anchor=west, right] at (1,-0.75) {(b)};
      \node[text width=1cm, anchor=west, right] at (1.25,-3) {(c)};
      
    \end{tikzpicture}
    
  \end{center}
  \caption{The Hasse diagram for a general $d$-dimensional SCS $\Gamma\subseteq\cP(\Sigma^{F_k^d})$, where (a) follows from \eqref{eq:ccaptcap}, (b) follows from Theorem \ref{indentropy1}, and (c) follows from \eqref{eq:hindohind}. }
  \label{fig:general}
\end{figure}
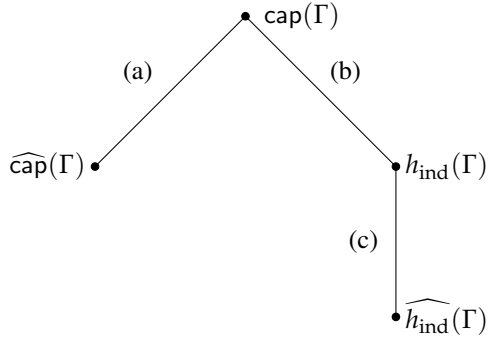

We note here that following the same arguments used in the proof of Theorem \ref{th:equ} would show that $\tcap(\God)\geq\tcap(\Gamma^{\otimes d+1})$ which means that in $\limsupup{d}\tcap(\God)$ the limit actually exists and equals to $\inf_d \tcap(\God)$.

We would also like to compare our results, as they apply to a specific
case study described in \cite{EliMeySch16}. Let
$\Gamma\subseteq\cP(\Sigma^k)$ be a convex one-dimensional SCS, and
recall that the axial product $\God$ is defined as
\[\God \eqdef \mathset{\mu\in\cP(\Sigma^{F_k^d}) ~:~ \forall i\in[d],\; \pi_{[k]\bunt_i}(\mu) \in\Gamma},\]
and thus
\[ \cB_n\parenv{\God}= \mathset{w\in F_n^d ~:~ \forall i\in[d],\; \fr^{[k]\bunt_i}_w\in\Gamma}.\]
The SCSs studied in \cite{EliMeySch16} were an averaged version of the axial product, namely,
\[\oGod\eqdef\mathset{\mu\in\cP(\Sigma^{F_k^d}) ~:~ \frac{1}{d}\sum_{i\in[d]}\pi_{[k]\bunt_i}(\mu) \in\Gamma},\]
and thus
\[ \cB_n\parenv{\oGod}= \mathset{w\in F_n^d ~:~ \frac{1}{d}\sum_{i\in[d]}\fr^{[k]\bunt_i}_w\in\Gamma}.\]
By convexity, it easily follows that
\[ \cB_n\parenv{\God} \subseteq \cB_n\parenv{\oGod},\]
and thus
\[ \tcap(\God) \leq \tcap(\oGod).\]

We now focus on the simple example known as the $(0,k,p)$-RLL SCS over
the binary alphabet $\Sigma=\mathset{0,1}$, which was the case study
of \cite{EliMeySch16}. The one-dimensional $(0,k,p)$-RLL SCS, $0\leq
p\leq 1$, is defined by
\begin{equation}
  \label{eq:rll}
  \Gamma_{k,p}\eqdef\mathset{\mu\in\cP(\Sigma^{k+1}) ~:~ \mu(1^{k+1})\leq p},
\end{equation}
where $1^{k+1}$ denotes the all-ones string of length $k+1$. This
example is a generalization of the well known inverted $(0,k)$-RLL
fully constrained system, since if we take $p=0$ we obtain the
inverted $(0,k)$-RLL. In \cite{EliMeySch16}, the authors found lower
and upper bounds on the internal capacity of $\oGodk$. We recall the relevant
lower bound here.

\begin{theorem}\cite[Th.~20]{EliMeySch16}
  \label{th:cbound}
  Let $\Gamma_{k,p}$ denote the one-dimensional $(0,k,p)$-RLL SCS given in
  \eqref{eq:rll}. Then, for all $0\leq p\leq \frac{1}{2^{k+1}}$,
  \[\ccap(\oGodk) \geq 1+d\parenv{\ccap(\Gamma)-1},\]
  whereas for all $\frac{1}{2^{k+1}}\leq p\leq 1$, $\ccap(\oGodk)=1$.
\end{theorem}

We first note that this theorem implies a lower bound on $\tcap(\oGodk)$,
\[  \tcap(\oGodk)\geq \ccap(\oGodk) \geq 1+d\parenv{\ccap(\Gamma)-1}.\]
The lower bound of \cite{EliMeySch16} eventually becomes negative, as
the dimension $d$ grows, and therefore, degenerate. However, using the
results of this paper,
\[ \tcap(\oGodk) \geq \hind(\Gamma_{k,p}),\]
and this bound does not depend on the dimension, and therefore, does
not degenerate. We provide an explicit numerical example:

\begin{example}
  Let us take $k=2$, and $p=0.05$, meaning that we restrict the
  frequency of the pattern $111$ to be at most $0.05$. Fix $d=3$. The
  lower bound on $\tcap(\oGodk)$ from \cite{EliMeySch16} uses
  $\ccap(\Gamma_{k,p})$. The latter can be calculated by solving an
  optimization problem using a computer. We obtain that
  $\ccap(\Gamma_{k,p})\approx 0.976$ which means that
  \[\tcap(\oGodk)\geq 1+3\cdot(0.976-1)\approx 0.928.\]
  Using the results of this paper, we use $\hind(\Gamma_{k,p})$ as a lower
  bound to $\tcap(\oGodk)$. Finding the supremum involved in the
  definition of $\hind(\Gamma_{k,p})$ is also not easy, and we lower bound it by
  guessing a specific measure.  We take each coordinate to be
  i.i.d.~Bernoulli $\sqrt[3]{0.05}$, and we get
  \[\tcap(\oGodk)\geq \hind(\Gamma_{k,p})\geq H_2(\sqrt[3]{0.05})\approx 0.949,\]
  which is a better lower bound than that of \cite{EliMeySch16}. Note
  that the upper bound gives $\ccap(\Gamma')\leq 0.983$. We further
  mention that the lower bound of \cite{EliMeySch16} gets increasingly
  worse as the dimension grows. For example, when $d=10$ we obtain by
  Theorem \ref{th:cbound} that $\tcap(\oGodk)\geq 0.76$ whereas using
  the independence entropy, the bound stays the same, i.e.,
  $\tcap(\oGod)\geq 0.949$. Finally, for all $d\geq 42$, the lower
  bound of \cite{EliMeySch16} becomes degenerate.
\end{example}

We present another example for $(0,1,p)$ with a more elaborate lower
bound.

\begin{figure}
 \centering
 \begin{subfigure}[t]{0.4\textwidth}
   \centering
   \begin{overpic}
     {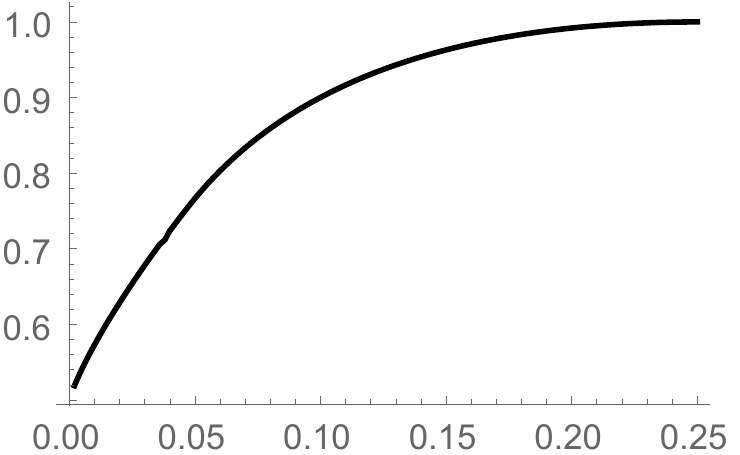}
     \put(100,7){$p$}
     \put(0,70){$\max_{x,y}\frac{1}{2n}H(\mu_{2n})$}
   \end{overpic}
   \caption{}
 \end{subfigure}%
 ~ 
 \begin{subfigure}[t]{0.6\textwidth}
   \centering
   \begin{overpic}
     {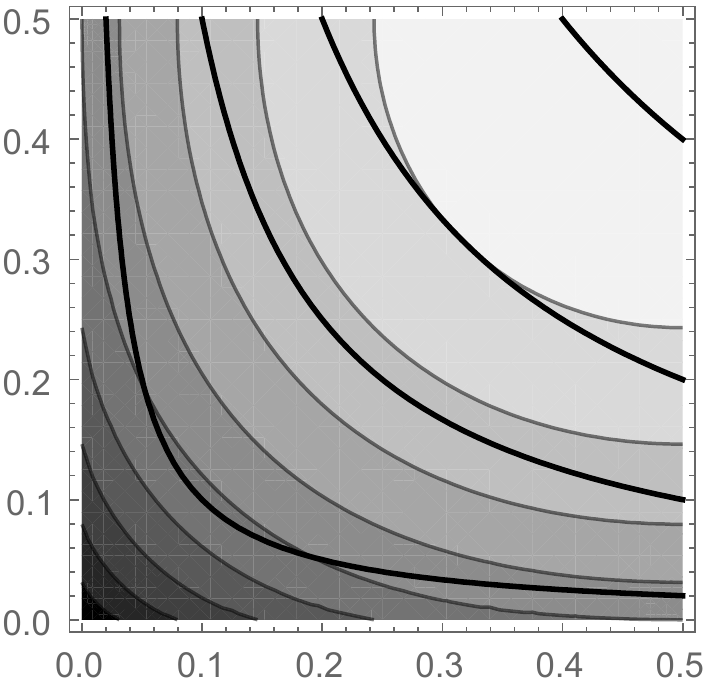}
     \put(100,7){$x$}
     \put(10,100){$y$}
     \put(100,12){$p=0.01$}
     \put(100,25){$p=0.05$}
     \put(100,42){$p=0.1$}
     \put(100,76){$p=0.2$}
   \end{overpic}
   \caption{}
 \end{subfigure}
 \caption{A lower bound on $\limsupup{d}\tcap(\Gamma_{1,p}^{\boxtimes d})$ is shown in (a), where
   (b) shows a contour plot of $\frac{1}{2}(H_2(x)+H_2(y))$ as well as the curves $xy=p$ for $p=0.01,0.05,0.1,0.2$.}
 \label{fig:rllp}
\end{figure}
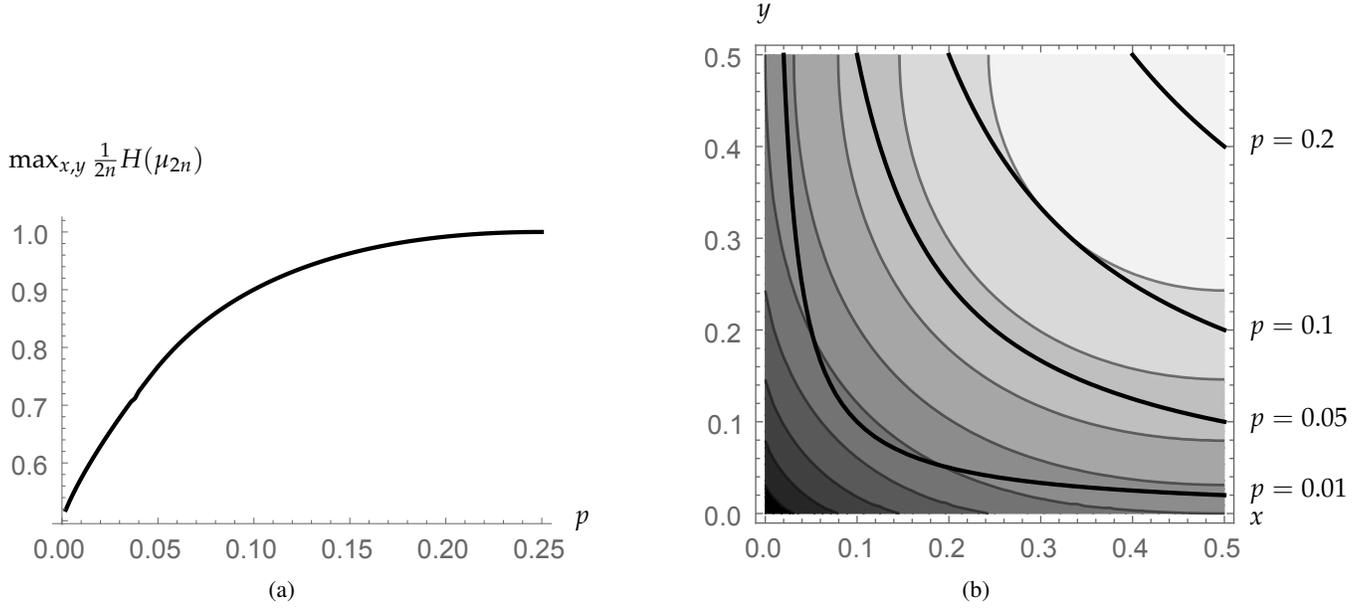

\begin{example}
  Take $k=1$ and consider $\Gamma_{k,p}$. From the results of this paper,
  \[ \limsupup{d}\tcap(\Gamma_{1,p}^{\boxtimes d})\geq \limsupup{d}\tcap(\God_{1,p})=\ohind(\Gamma_{1,p})\geq \hind(\Gamma_{1,p}).\]
  We lower bound $\hind(\Gamma_{1,p})$ by devising a product measure
  $\mu_{2n}\in(\cP(\Sigma))^{2n}$, for all $n\in\N$. The measures use
  two parameters $0\leq x,y\leq 1$, using a $\mathrm{Bernoulli}(x)$
  distribution for positions with odd indices, and a
  $\mathrm{Bernoulli}(y)$ for positions with even indices. Thus,
  \[ \hind(\Gamma_{1,p})\geq \max_{x,y}\frac{1}{2n}H(\mu_{2n})=\max\mathset{\frac{1}{2}(H_2(x)+H_2(y)) ~:~ 0\leq x,y\leq 1,\ xy\leq p}.\]
  Due to monotonicity, the maximization problem always has a solution
  on the curve $xy=p$, which in the high range is unique
  $x=y=\sqrt{p}$, and in the lower range has two symmetric
  solutions. For example, for $p=0.2$ the optimal solution is
  $x=y=\sqrt{0.2}$. However, for $p=0.01$, the first optimal solution
  is $x\approx 0.454$, $y\approx 0.022$, and the symmetric solution is
  $x\approx 0.022$, $y\approx 0.454$. This is depicted in Figure
  \ref{fig:rllp}.

  We note that this bound agrees with the solution for the fully
  constrained case, $\limsupup{d}\tcap(\Gamma_{k,0}^{\boxtimes
      d})=\frac{1}{2}$ which was solved in \cite{MeyPav14}. We
  conjecture that Figure \ref{fig:rllp}(a) indeed shows the exact
  limiting capacity.
\end{example}

\appendices
\section{Cyclic and Non-cyclic Capacities}
\label{appA} 

The goal of this appendix is to show that the capacity, as we defined
it cyclically, equals the (traditionally non-cyclic) capacity in the
case of fully constrained systems. 
	
\begin{definition}
  \label{def:fullyold}
  Let $d,k\in\N$. A \emph{(traditional) fully constrained system} is a set
  $\Phi\subseteq \Sigma^{F_k^d}$ of $d$-dimensional words, called
  \emph{forbidden patterns}. The set of all \emph{admissible} words in
  $\Sigma^{F_n^d}$ is defined as
  \[\cBcom_n(\Phi) \eqdef \mathset{ x\in\Sigma^{F_n^d} ~:~ \forall \bv\in
    F_{n-k}^d,\; x_{\bv+F_k^d} \notin \Phi}.\]
  The (combinatorial) \emph{capacity} of $\Phi$ is defined by
  \[\capcom(\Phi)\eqdef \limsupup{n}\frac{1}{|F_n^d|}\log_2 \abs{\cBcom_n(\Phi)}.\]
\end{definition}

Intuitively, a traditional fully constrained system is a set of words that
do not contain any forbidden pattern \emph{non-cyclically}. Given a
(traditional) fully constrained system $\Phi\subseteq\Sigma^{F_k^d}$, we
can construct a set of measures $\Gamma_\Phi$ defined as follows,
\begin{equation}
  \label{eq:gamphi}
  \Gamma_\Phi \eqdef \mathset{ \mu\in\cP(\Sigma^{F_k^d}) ~:~ \mu(\Sigma^{F_k^d}\setminus\Phi)=1}.
\end{equation}
Thus, $\Gamma_\Phi$ is a SCS which is fully constrained in the sense
of Definition \ref{def:fullynew}. Since Definition \ref{def:fullynew}
is more restrictive, by requiring forbidden patterns to not appear in
admissible words \emph{cyclically}, we immediately have
\[ \cB_n(\Gamma_\Phi) \subseteq \cBcom_n(\Phi),\]
implying also
\[ \ccap(\Gamma_\Phi) \leq \capcom(\Phi).\]
However, we now prove that the capacity of $\Gamma_\Phi$ does
equal the (combinatorial) capacity of $\Phi$.

\begin{prop}
  Let $d,k\in\N$. Let $\Phi\subseteq \Sigma^{F_k^d}$ be a fully
  constrained system as in Definition \ref{def:fullyold}, and let
  $\Gamma_\Phi\subseteq\cP(\Sigma^{F_k^d})$ be its corresponding fully
  constrained system as in Definition \ref{def:fullynew}.  If
  $\cBcom_n(\Phi)\neq\emptyset$ for all large enough $n\in\N$,
  then
  \[\tcap(\Gamma_\Phi)=\capcom(\Phi).\]
\end{prop}
		
\begin{IEEEproof}
  We first show that $\capcom(\Phi)\leq \tcap(\Gamma_\Phi)$. Fix
  $\epsilon>0$, and for $n\in\N$, $n\geq k$, consider the $k$-boundary
  of $F_n^d$ which is defined as $F_n^d\setminus F_{n-k}^d$. Note that
  $|F_n^d\setminus F_{n-k}^d|=n^d-(n-k)^d$. Let
  $w\in\cBcom_n(\Phi)$. While $w$ does not contain any forbidden pattern when considering the coordinates non-cyclically, it may contain some when considering the coordinates cyclically. The number
  of occurrences of forbidden patterns (cyclically) in $w$ is at most
  $|F_n^d\setminus F_{n-k}^d|=n^d-(n-k)^d$. For all large enough $n$ we have
  $\frac{n^d-(n-k)^d}{n^d}\leq\epsilon$, hence
  \[\cBcom_n(\Phi)\subseteq \cB_n(\Be(\Gamma_\Phi)).\]
  Thus, for every
  $\epsilon>0$,
  \[\capcom(\Phi)\leq \ccap(\Be(\Gamma_\Phi)).\] 
  Taking the limit as $\epsilon\to 0$ we obtain 
  \[\capcom(\Phi)\leq \tcap(\Gamma_\Phi).\]
    
  In the other direction, we now show that
  $\tcap(\Gamma_\Phi)\leq \capcom(\Phi)$. Let $\delta_0>0$ and take
  $n_0\in\N$ large enough such that
  \[\frac{1}{n_0^d}\log_2 \abs{\cBcom_{n_0}(\Phi)}\leq \capcom(\Phi)+\frac{1}{3}\delta_0.\]
  Denote the number of forbidden patterns by $t\eqdef\abs{\Phi}$.
  Take $\delta>0$ small enough such that both
  \[\frac{t(1+\delta)}{n_0^d}H_2\parenv{\frac{\delta}{1+\delta}}\leq \frac{1}{3}\delta_0, \qquad\text{and}\qquad t\delta\log_2\abs{\Sigma}\leq \frac{1}{3}\delta_0, \]
  where $H_2(\cdot)$ is the binary entropy function. Finally, for
  every $n\geq n_0$, denote $m\eqdef \floorenv{n/n_0}$, and choose any
  $0 < \epsilon \leq \delta/n_0^d$.

  Consider a word $w\in\cB_{n}(\Be(\Gamma_\Phi))$. We say $w$ is made up a
  concatenation of $m^d$ $F_{n_0}^d$-blocks, namely, a block is a set
  of positions $n_0\bv + F_{n_0}^d$, where $\bv\in F_m^d$, as well a
  boundary, namely, the set of positions $F_n^d\setminus F_{mn_0}^d$.
  By our choice of parameters, the number of occurrences (perhaps
  cyclically) of any forbidden pattern from $\Phi$ is at most
  \[\epsilon|F_n^d| \leq \epsilon (m+1)^d n_0^d \leq \delta (m+1)^d.\]
  This serves also as an upper bound on the number of blocks fully
  containing (non-cyclically) this forbidden pattern. Since there are
  $t$ forbidden patterns, the number of blocks that are devoid
  (non-cyclically) of any forbidden pattern, is at least $m^d -
  t\delta(m+1)^d$. Such blocks are in fact words from
  $\cBcom_{n_0}(\Phi)$.
  
  Fixing a specific type of forbidden pattern, and considering each
  occurrence of it as a ball, we have at most $\delta(m+1)^d$ balls,
  which we throw into $m^d+1$ bins ($m^d$ blocks, and another
  ``virtual'' bin for patterns that are not fully contained within a
  single block).  The total number of ways to throw these ball into
  bins is exactly
  $\binom{m^d+1+\delta(m+1)^d}{\delta(m+1)^d}$. Raising this to the
  power of $t$ gives an upper bound on the number of ways the $t$
  forbidden patterns are dispersed among the blocks. In total we have,
  \begin{align*}
    \abs{\cB_{n}(\Be(\Gamma_\Phi))} &\leq \binom{m^d+1+\delta(m+1)^d}{\delta(m+1)^d}^t  \abs{\cBcom_{n_0}(\Phi)}^{m^d-t(m+1)^d\delta} \abs{\Sigma}^{t\delta(m+1)^d n_0^d} \abs{\Sigma}^{n^d-(mn_0)^d},
  \end{align*} 
  where the binomial coefficient follows from upper bounding the way
  forbidden patterns are dispersed among blocks, the following term
  counts the number of ways to fill blocks that do not contain
  (non-cyclically) any forbidden word, and the last term counts the
  ways to arbitrarily fill in the rest of the positions. Thus,
  \begin{align*}
    \ccap(\Be(\Gamma_\Phi))&=\limsupup{n}\frac{1}{n^d}\log_2 \abs{\cB_{n}(\Be(\Gamma_\Phi))} \\
    &\leq
  \frac{t(1+\delta)}{n_0^d}H_2\parenv{\frac{\delta}{1+\delta}}
  +\frac{1}{n_0^d}\log_2 \abs{\cBcom_{n_0}(\Phi)} + t\delta\log_2\abs{\Sigma}\\
  &\leq \delta_0 + \capcom(\Phi).
  \end{align*}
  Taking the limit as $\epsilon\to 0$, we get
  \[ \tcap(\Gamma_\Phi) \leq \delta_0 + \capcom(\Phi).\]
  Finally, since this holds for any $\delta_0>0$, we get the desired result,
  \[ \tcap(\Gamma_\Phi) \leq \capcom(\Phi).\]
\end{IEEEproof}

\section{Independence Entropy for Fully Constrained Systems}
\label{appB}

Here we Prove Theorem \ref{indentropyequivalent}.
We begin by recalling relevant definitions from \cite{LouMarPav13}.  A
$\Z^d$ shift space $X$, is a subset $X\subseteq \Sigma^{\ZD}$ that is
closed under shifts, i.e., for all $\bv\in\ZD$, and all $x\in X$,
$\sigma_{\bv}(x)\in X$.

\begin{definition}
  \label{def:fullyshift}
  Let $d,k\in\N$. Given a set of forbidden words $\Phi\subseteq
  \Sigma^{F_k^d}$, the $\Z^d$ shift space over $\Sigma$ defined by
  $\Phi$ is
  \[ X_\Phi \eqdef \mathset{ x\in\Sigma^{\Z^d} ~:~ \forall \bv\in
    \Z^d, x_{\bv+F_k^d} \notin \Phi}.\]
\end{definition}

Given a finite alphabet $\Sigma$, let $\tilde{\Sigma}$ denote the set of
all non-empty subset of $\Sigma$, i.e.,
\[ \tilde{\Sigma} \eqdef \mathset{ A\subseteq \Sigma ~:~ A\neq\emptyset}.\]

\begin{definition}
  Let $d\in\N$, $S\subseteq \ZD$, and let $\tilde{x}$ be a configuration
  on $S$ over $\tilde{\Sigma}$, i.e.,
  $\tilde{x}\in\tilde{\Sigma}^{S}$. Denote by $\varphi(\tilde{x})$ the set of
  fillings of $\tilde{x}$,
  \[\varphi(\tilde{x})\eqdef\mathset{x\in\Sigma^S ~:~ \forall \bv\in S,\; x_{\{\bv\}}\in\tilde{x}_{\{\bv\}}}.\]
\end{definition}

\begin{definition}
Let $d\in\N$, and let $X$ be a $\ZD$ shift space over $\Sigma$. We
denote by $\tilde{X}$ the multi-choice shift space corresponding to $X$,
\[\tilde{X}\eqdef\mathset{\tilde{x}\in\tilde{\Sigma}^{\ZD} ~:~ \varphi(\tilde{x})\subseteq X}.\]
We also denote by $\cB_n(\tilde{X})$ the set of all eligible
configurations on $F_n^d$ in $\tilde{X}$, i.e.,
\[\cB_n(\tilde{X}) \eqdef \mathset{ \tilde{x}_{F_n^d} ~:~ \tilde{x}\in\tilde{X} }.\]
\end{definition}

\begin{definition}
  Let $d\in\N$, and let $X$ be a $\ZD$ shift space. We define the
  \emph{combinatorial independence entropy of $X$}, denoted as
  $\hindcom(X)$, by
  \[\hindcom(X)\eqdef\limsupup{n} \frac{1}{n^d}\max \mathset{\log_2 \abs{\varphi(\tilde{w})} ~:~ \tilde{w}\in\cB_n(\tilde{X}) }.\]
\end{definition}

Note that in \cite{LouMarPav13} the definition of combinatorial
independence entropy is slightly more general and defined over all
shapes and not only on the shapes $F_n^d$. Finally, given a fully
constrained system $\Phi\subseteq \Sigma^{F_k^n}$ (see Definition
\ref{def:fullyold}), its representation as a SCS is given by
$\Gamma_\Phi$ in \eqref{eq:gamphi}. We are now ready to prove Theorem
\ref{indentropyequivalent}.

\begin{IEEEproof}[Proof of Theorem \ref{indentropyequivalent}]
  Let $d,k\in\N$, and let $\Phi\subseteq \Sigma^{F_k^d}$ be a fully
  constrained system, with its SCS representation $\Gamma_\Phi$ from
  \eqref{eq:gamphi}. The claim we want to prove is that
  \[ \hindcom(X_\Phi) = \ohind(\Gamma_\Phi).\]
  
  First, we show that $\hindcom(X_{\Phi})\leq\ohind(\Gamma_\Phi)$. For
  every $n\in\N$ choose $\tilde{w}_n\in\cB_n(\tilde{X}_{\Phi})$ which
  maximizes $\abs{\varphi(\tilde{w}_n)}$. Now consider the independent
  measures $\mu_n$ such that $\pi_{\mathset{\bv}}(\mu_n)$ is the
  uniform distribution over $(\tilde{w}_n)_{\mathset{\bv}}$. Note that
  in $\cB_n(\tilde{X})$, the forbidden patterns are considered without
  modulo while in $\cB_n(\Gamma_\Phi)$ the calculation of the
  marginals' average uses modulo $n$. Therefore, if a filling
  $\varphi(\tilde{w}_n)$ belongs to $X_{\Phi}$, in $\mu_n$ there is
  perhaps a positive probability to see a forbidden pattern only in
  the boundaries. In $F_n^d$, the $k$-boundary is the set
  $F_n^d\setminus F_{n-k}^d$ of size $n^d-(n-k)^d$. Since
  $(n^d-(n-k)^d)/n^d \to 0$ as $n\to\infty$, we obtain that for
  every $\epsilon>0$, for every $n\in\N$ such that
  $(n^d-(n-k)^d)/n^d\leq \epsilon$, we have that
  $\mu_n\in\Be(\Gamma_\Phi)$. Thus,
  \begin{align*}
  	\hind(\Be(\Gamma_\Phi)) &= \limsupup{n} \sup_{\mu\in\ocP_n(\Be(\Gamma_\Phi))} \frac{1}{n^d}H(\mu)\\
  	& \geq \limsupup{n}\frac{1}{n^d}H(\mu_n)\\
  	& = \limsupup{n}\frac{1}{n^d}\log_2\abs{\varphi(\tilde{w}_n)}\\
  	& = \hindcom(X_\Phi).
  \end{align*}
  Taking $\epsilon\to 0$ we obtain 
  \[\ohind(\Gamma_\Phi)\geq \hindcom(X_\Phi).\]
  
  We now show that $\ohind(\Gamma_\Phi)\leq \hindcom(X_\Phi)$. Fix
  $\delta>0$ and take $\delta_1>0$ small enough such that
  $\delta_1<\frac{1}{3}\delta$. Take $n_0\in\N$ large enough such that
  for all $n\geq n_0$,
  \begin{equation}
    \label{eq:closeenough}
  \frac{1}{n^d}\max_{\tilde{w}\in\cB_{n}(\tilde{X}_{\Phi})}
  \mathset{\log_2 \abs{\varphi(\tilde{w})}}\leq\hindcom(X_\Phi)+\frac{1}{3}\delta.
  \end{equation}
  We now take $\epsilon>0$ small enough such that all the following hold,
  \begin{align}
    -\abs{\Sigma}\sqrt[k^d]{n_0^d \epsilon^{\frac{1}{4}}}\log_2 \sqrt[k^d]{n_0^d \epsilon^{\frac{1}{4}}}&<\frac{1}{3}\delta, \label{eq:entdec} \\
    2^d\epsilon^{\frac{3}{4}}\log_2 \abs{\Sigma} &<\frac{1}{2}\delta_1, \label{eq:2eps}\\
    \abs{\hind\parenv{\Be(\Gamma_\Phi)}-\ohind(\Gamma_\Phi)} &\leq \frac{1}{16}\delta_1. \nonumber
  \end{align}
  By the definition of $\hind\parenv{\Be(\Gamma_\Phi)}$ we may find
  $n\geq n_0$ large enough such that all the following hold,
  \begin{align}
    2\parenv{1-\parenv{1-\frac{n_0}{n}}^d}\log_2\abs{\Sigma} &\leq \frac{1}{4}\delta_1, \label{eq:nchoice}\\
    \abs{\sup_{\mu\in \ocP_n(\Be(\Gamma_\Phi))}\frac{1}{n^d}H(\mu)-\ohind(\Gamma_\Phi)} &\leq \frac{1}{8}\delta_1, \nonumber
  \end{align}
  and there exists $\mu\in\ocP_n(\Be(\Gamma_\Phi))$ for which 
  \[\left|\frac{1}{n^d}H(\mu)-\ohind(\Gamma_\Phi)\right|\leq \frac{1}{4}\delta_1.\]
  Since $\mu\in\ocP_n(\Be(\Gamma_\Phi))$, we have
  \[\frac{1}{n^d}\sum_{\bv\in F_n^d}\pi_{F_k^d+\bv}(\mu)(\Phi)\leq \epsilon.\]
  
  Denote $m\eqdef \floorenv{n/n_0}$. We now partition $F_n^d$ into
  $m^d$ blocks of shape $F_{n_0}^d$ in the natural way,
  $\mathset{n_0\bv+F_{n_0}^d : \bv\in F_m^d}$, as well as a boundary
  $F_n^d\setminus F_{mn_0}^d$. Note that
  \[\mu\cong \bigotimes_{\bv\in
    F_m^d}\pi_{n_0\bv+F_{n_0}^d}(\mu) \otimes\pi_{F_n^d\setminus
    F_{mn_0}^d}(\mu).
  \]
  Since $\mu$ is independent we obtain
  \begin{align}
    \abs{ \frac{1}{(mn_0)^d}H(\pi_{F_{mn_0}^d}(\mu))- \ohind(\Gamma_\Phi)} &\leq 
    \abs{ \frac{1}{(mn_0)^d}H(\pi_{F_{mn_0}^d}(\mu))- \frac{1}{n^d}H(\mu)} +\abs{ \frac{1}{n^d}H(\mu)-\ohind(\Gamma_\Phi)} \nonumber \\
  	&\leq (mn_0)^d\parenv{\frac{1}{(mn_0)^d}-\frac{1}{n^d}} \log_2 \abs{\Sigma}+\frac{n^d-(mn_0)^d}{n^d}\log_2 \abs{\Sigma}+ \frac{1}{4}\delta_1 \nonumber \\
    &\leq 2\frac{n^d-(mn_0)^d}{n^d}\log_2 \abs{\Sigma}+ \frac{1}{4}\delta_1 \nonumber \\
    &\leq 2\parenv{1-\parenv{1-\frac{n_0}{n}}^d}\log_2\abs{\Sigma}+ \frac{1}{4}\delta_1 \nonumber \\
    &\leq \frac{1}{2}\delta_1, \label{eq:halfohind}
  \end{align}
  where the last inequality holds due to \eqref{eq:nchoice}. Let
  $Z:F_m^d\to\R$ be a function defined by
  \[Z(\bv)\eqdef\frac{1}{n_0^d} \sum_{\mathbf{u}\in F_{n_0}^d}\pi_{F_k^d+n_0\bv+\mathbf{u}}(\mu)(\Phi)\]
  (with coordinates taken modulo $n$). Note that since
  $\mu\in\ocP_n(\Be(\Gamma_\Phi))$, we have
  \[\frac{1}{n^d}\sum_{\bv\in F_n^d}\pi_{F_k^d+\bv}(\mu)(\Phi)\leq \epsilon.\]
  If we now take $\bv$ to be random uniformly distributed in $F_m^d$, then
  \[E[Z(\bv)]= \frac{1}{m^d}\sum_{\bv\in F_m^d} \frac{1}{n_0^d} \sum_{\mathbf{u}\in F_{n_0}^d}\pi_{F_k^d+n_0\bv+\mathbf{u}}(\mu)(\Phi) \leq \parenv{1+\frac{1}{m}}^d\epsilon.\]
  By Markov's inequality we have
  \begin{equation}
    \label{eq:mark}
    \Pr\parenv{Z(\bv)\geq\epsilon^{\frac{1}{4}}} \leq \epsilon^{-\frac{1}{4}}E[Z(\bv)]\leq \parenv{1+\frac{1}{m}}^d\epsilon^{\frac{3}{4}}.
  \end{equation}

  Recall that each $\bv\in F_m^d$ may be identified with the
  $F_{n_0}^d$ block of $F_n^d$ in coordinates $n_0\bv+F_{n_0}^d$. Define,
  \[ \cL \eqdef \mathset{\bv\in F_m^d ~:~ Z(\bv)\geq \epsilon^{\frac{1}{4}}}.\]
  Since $\bv$ was distributed uniformly in $F_m^d$, by \eqref{eq:mark} we have,
  \begin{equation}
    \label{eq:cl}
    \abs{\cL} \leq (m+1)^d\epsilon^{\frac{3}{4}}.
  \end{equation}
  It now follows that
  \begin{align*}
    \frac{1}{(mn_0)^d}\sum_{\bv\in F_m^d\setminus\cL} H(\pi_{n_0\bv+F_{n_0}^d}(\mu)) &=
    \frac{1}{(mn_0)^d}H(\pi_{F_{mn_0}^d}(\mu)) - \frac{1}{(mn_0)^d}\sum_{\bv\in \cL} H(\pi_{n_0\bv+F_{n_0}^d}(\mu))\\
    &\overset{(a)}{\geq} \ohind(\Gamma_\Phi)-\frac{1}{2}\delta_1 - \frac{1}{(mn_0)^d}\sum_{\bv\in \cL} H(\pi_{n_0\bv+F_{n_0}^d}(\mu))\\
    &\overset{(b)}{\geq} \ohind(\Gamma_\Phi)-\frac{1}{2}\delta_1 - \frac{(m+1)^d\epsilon^{\frac{3}{4}}n_0^d}{(mn_0)^d}\log_2\abs{\Sigma}\\
    &\geq \ohind(\Gamma_\Phi)-\frac{1}{2}\delta_1 - 2^d\epsilon^{\frac{3}{4}}\log_2\abs{\Sigma}\\
    &\overset{(c)}{>} \ohind(\Gamma_\Phi)-\delta_1 \\
    &> \ohind(\Gamma_\Phi)-\frac{1}{3}\delta,
  \end{align*}
  where $(a)$ follows from \eqref{eq:halfohind}, $(b)$ follows from
  \eqref{eq:cl}, and $(c)$ follows from \eqref{eq:2eps}. Since there
  are at most $m^d$ summands on the left-hand side, there exists
  $\bv_0\in F_m^d\setminus \cL$ such that
  \begin{equation}
    \label{eq:ohindmd3}
  \frac{1}{n_0^d} H(\pi_{n_0\bv_0+F_{n_0}^d}(\mu)) \geq \ohind(\Gamma_\Phi)-\frac{1}{3}\delta.
  \end{equation}
  We denote by $\nu$ the independent measure
  $\nu\eqdef\pi_{F_{n_0}^d+n_0\bv_0}(\mu)$.

  Note that if we consider $\nu$ in a non-cyclic manner, we obtain
  that
  \[\frac{1}{(n_0-k+1)^d}\sum_{\mathbf{u}\in
    F_{n_0-k+1}^d}\pi_{\bu+F_k^d}(\nu)(\Phi)\leq
  \frac{n_0^d}{(n_0-k+1)^d}\epsilon^{\frac{1}{4}},
  \]
  and in particular, for every coordinate $\mathbf{u}\in
  F_{n_0-k+1}^d$, we have that $\pi_{\bu+F_k^d}(\nu)(\Phi)\leq
  n_0^d\epsilon^{\frac{1}{4}}$. Let us define
  \[p\eqdef \sqrt[k^d]{n_0^d\epsilon^{\frac{1}{4}}}.
  \]
  Hence, since $\nu$ is an independent measure, if $a\in\Phi$ then
  there must be a coordinate $\mathbf{t}\in F_k^d$ for which
  $\pi_{\mathbf{u+t}}(\nu)(a_{\mathbf{t}})\leq p$.

  We now construct a configuration $\tilde{w}\in
  \tilde{\Sigma}^{F_{n_0}^d}$. For every coordinate $\mathbf{u}\in
  F_{n_0}^d$ we take
  \[ \tilde{w}_{\bu}= \mathset{ a\in\Sigma ~:~ \pi_{\mathset{\bu}}(\nu)(a)> p}.\]
  By our previous observation,
  $\tilde{w}\in\cB_{n_0}(\tilde{X}_{\Phi})$ since any filling of
  $\tilde{w}$ cannot contain a forbidden word from $\Phi$ as it
  requires at least one position $\bu$ such that
  $\pi_{\mathset{\bu}}(\nu)\leq p$. Moreover,
  \begin{align*}
    \log_2\abs{\tilde{w}_{\bu}} &\geq -\sum_{a\in\tilde{w}_{\bu}} \pi_{\mathset{\bu}}(\nu)(a)\log_2 (\pi_{\mathset{\bu}}(\nu)(a))\\
    & = H(\pi_{\mathset{\bu}}(\nu)) + \sum_{a\in\Sigma\setminus \tilde{w}_{\bu}}\pi_{\mathset{\bu}}(\nu)(a)\log_2 (\pi_{\mathset{\bu}}(\nu)(a)) \\
    &\geq H(\pi_{\mathset{\bu}}(\nu)) + \abs{\Sigma}p\log_2 p\\
    &> H(\pi_{\mathset{\bu}}(\nu)) -\frac{1}{3}\delta
  \end{align*}
  where the last inequality follows from \eqref{eq:entdec}. Hence, using \eqref{eq:ohindmd3},
  \[ \frac{1}{n_0^d}\log_2\abs{\varphi(\tilde{w})} = \frac{1}{n_0^d}\sum_{\bu\in F_{n_0}^d} \log_2\abs{\tilde{w}_{\bu}}=\frac{1}{n_0^d}\sum_{\bu\in F_{n_0}^d} H(\pi_{\mathset{\bu}}(\nu)) -\frac{1}{3}\delta \geq \ohind(\Gamma_\Phi)-\frac{2}{3}\delta. \]
  Finally, using \eqref{eq:closeenough}, this implies that,
  \[\hindcom(X_\Phi) \geq \frac{1}{n_0^d}\max_{\tilde{w}\in\cB_{n_0}(\tilde{X}_{\Phi})} \mathset{\log_2 \abs{\varphi(\tilde{w})}} - \frac{1}{3}\delta \geq \ohind(\Gamma_{\Phi})-\delta.\]
  Since this holds for every $\delta>0$ we have $\ohind(\Gamma_\Phi)\leq \hindcom(X_\Phi)$, as claimed.
\end{IEEEproof}

\bibliographystyle{IEEEtranS}
\bibliography{allbib}

\end{document}